\documentclass[12pt]{article}

\global\arraycolsep=1pt
\usepackage{amsmath}
\usepackage{amssymb}
\usepackage {lscape}
\newcommand{\appSymm}{A}
\newcommand{\appProofHL}{B}
\newcommand{\appMacOp}{C}
\newcommand{\appProofMO}{D}
\newcommand{\appTorus}{E}
\newcommand{\appPositivity}{F}
\newcommand{\be}{\begin{equation}}
\newcommand{\ee}{\end{equation}}
\newcommand{\ba}{\begin{eqnarray}}
\newcommand{\ea}{\end{eqnarray}}
\newcommand{\ds}{\displaystyle}
\renewcommand{\theequation}{\thesection.\arabic{equation}}
\newcommand{\Section}{\setcounter{equation}{0} \section}
\newcommand{\Appendix}[2]{
\renewcommand{\theequation}{#1.\arabic{equation}}\setcounter{equation}{0}
\renewcommand{\thesubsection}{#1.\arabic{subsection}}\setcounter{subsection}{0}
\section*{Appendix #1. #2}
}
\renewcommand{\thefootnote}{\fnsymbol{footnote}}
\newcommand{\bC}{{\mathbb C}}
\newcommand{\bN}{{\mathbb N}}
\newcommand{\bZ}{{\mathbb Z}}
\newcommand{\cS}{{\cal S}}
\newcommand{\proof}{\noindent{\it Proof.\hskip10pt}} 
\newcommand{\qed}{\hfill\fbox{}}
\newcommand{\rhs}{{\rm rhs}}
\newcommand{\lhs}{{\rm lhs}}
\newcommand{\ha}{{1\over2}}
\newcommand{\Exp}[1]{\exp\left\{#1\right\}}

\newcommand{\qint}[2]{[\, #1\, ]_{#2}}
\newcommand{\qintf}[2]{[\, #1\, ]!_{#2}}
\newcommand{\qbin}[3]{\left[{#1\atop #2}\right]_{#3}}

\newcommand{\bfa}{\mathbf a}
\newcommand{\bq}{\mathbf q}
\newcommand{\bt}{\mathbf t}
\newcommand{\tildeq}{\tilde q}
\newcommand{\tildet}{\tilde t}
\newcommand{\tZ}{\widetilde Z}
\newcommand{\tU}{\widetilde U}
\newcommand{\tV}{\widetilde V}
\newcommand{\Diff}[1]{T_{q,x_{#1}}}
\newcommand{\vertex}{\varphi}
\newcommand{\z}{z}
\newcommand{\w}{w}
\newcommand{\tw}{\tilde w}
\newcommand{\Ii}[1]{I_{#1}}
\newcommand{\Ji}[1]{J_{#1}}
\newcommand{\txx}[2]{[\, #1,#2\, ]}
\newcommand{\Dp}{\Delta}
\newcommand{\HopfLink}[5]{\bar {\cal P}_{#1,#2}(#3;#4,#5)}
\newcommand{\tz}{z}
\newcommand{\pz}{{\mbox{\boldmath $z$}}}
\newcommand{\GIKVQ}{Q^{{\rm GIKV}}}
\newcommand{\viQ}{v^{-1}Q}
\newcommand{\Q} {Q}
\newcommand{\Qv}{Q}
\newcommand{\vQ}{Q}
\newcommand{\vitN}{t^N}
\newcommand{\tNqt}{}
\newcommand{\gqt}[1]{g_{#1}}
\newcommand{\gqq}[1]{f_{#1}}
\newcommand{\pos}[1]{}
\newcommand{\nega}[1]{#1}

\newcommand\FigContour{%
\begin{figure}
\begin{center}
\begin{picture}(120,30)(-60,-15)
\put(-50,0){$\times$}
\put(-35,0){$\times$}
\put(-25,0){$\cdots$}
\put(-15,0){$\times$}
\put( 0, 0){$\times$}
\put(15,0){$\times$}
\put(25,0){$\cdots$}
\put(35,0){$\times$}
\put(50,0){$\times$}
\put(-50,-5){$\infty$}
\put(-35,-5){${1\over x_N}$}
\put(-15,-5){${1\over x_1}$}
\put(0,-5){$z_r$}
\put(15,-5){${z_{r-1}\over t}$}
\put(35,-5){${z_1\over t}$}
\put(50,-5){$0$}
\put(-30,-1.5){\oval(60,27)[l]}
\put(-30, 5){\oval(63,14)[tr]}
\put(-30,-8){\oval(63,14)[br]}
\put(-28, 12){\vector(1,0){0}}
\put(-60,-1){\vector(0,1){0}}
\put(-30,-15){\vector(-1,0){0}}
\put(1.5,4){\vector( 0,-1){0}}
\end{picture}
\end{center}
\caption{The constant term in $z_r$ can be represented as the contour integral.}
\end{figure}
}

\begin{document}

%
\begin{titlepage}
\begin{flushright}
{Oct. 1, 2009}
\end{flushright}
\vspace{0.5cm}
\begin{center}
{\Large \bf
Macdonald operators \\ and \\
homological invariants of the colored Hopf link}
\vskip1.0cm
{\large Hidetoshi Awata and Hiroaki Kanno}
\vskip 1.0em
{\it 
Graduate School of Mathematics \\
Nagoya University, Nagoya, 464-8602, Japan}
\end{center}
\vskip1.5cm



\begin{abstract}
Using a power sum (boson) realization for the Macdonald operators,
we investigate the Gukov, Iqbal, Koz\c{c}az and Vafa (GIKV) proposal 
for the homological invariants of the colored Hopf link,
which include Khovanov-Rozansky homology as a special case.
We prove the polynomiality 
of the invariants obtained by GIKV's proposal for arbitrary representations.
We derive a closed formula of the invariants of the colored Hopf link 
for antisymmetric representations.
We argue that a little amendment of GIKV's proposal is required 
to make all the coefficients of the polynomial non-negative integers.
\end{abstract}
\end{titlepage}


\renewcommand{\thefootnote}{\arabic{footnote}} \setcounter{footnote}{0}



\Section{Introduction and notation} 



In the setup of string theory, the invariants of the colored Hopf link are
identified as topological open string amplitudes on the deformed
conifold $T^*S^3$ as follows: 
the $U(N)$ Chern-Simons theory is realized
by topological string on $T^*S^3$ with $N$ topological $D$-branes wrapping
on the base Lagrangian submanifold $S^3$.
The Hopf link in $S^3$ consisting of two knots ${\cal K}_1$ and
${\cal K}_2$ can be introduced by a pair of new $D$-branes
wrapping on Lagrangian three-cycles ${\cal L}_1$ and ${\cal L}_2$
such that $S^3 \cap {\cal L}_i = {\cal K}_i$ \cite{OV1}.
The topological open string amplitude
of this brane system is supposed to give the invariants of the
Hopf link. The coloring or the representation attached to
each knot ${\cal K}_i$ is related to the boundary states
of the open string ending on ${\cal L}_i$ by the Frobenius
relation. After geometric transition or by the large $N$ duality
\cite{GV, Va, OV2},
this brane configuration is mapped to the resolved
conifold ${\cal O}(-1) \oplus {\cal O}(-1) \to {\bf P}^1$.
The $D$-branes wrapping on $S^3$ disappear,
but a pair of Lagrangian $D$-branes remains as a remnant of the Hopf link.
We can describe the resulting $D$-brane system in terms of the toric diagram
and compute the corresponding amplitude by the method of a topological vertex
 \cite{AMV, AKMV}.


The aim of this paper is to investigate 
the conjecture of Gukov, Iqbal, Koz\c{c}az and Vafa (GIKV) \cite{rf:GIKV} 
on the superpolynomial $\overline{\cal P}_{\lambda, \mu} (\bfa, \bq,\bt)$ 
of the homological invariants of the Hopf link colored by two representations
$\lambda$ and $\mu$.
The superpolynomial of our interest 
is a polynomial in $(\bfa,\bq,\bt)\in\bC^3$, 
such that a specialization $\bfa = \bq^N$ leads 
the Poincar\'e polynomial of  the $\mathfrak{sl} (N)$ link homology 
$\overline{\cal P}_{\mathfrak{sl}(N); \lambda, \mu} (\bq,\bt) $, 
which is a two parameter $(\bq,\bt)$ version of the $\mathfrak{sl}(N)$ link invariants.
When the coloring is the $N$ dimensional 
defining representation, it is called the Khovanov-Rozansky homology \cite{rf:KR}. 
In \cite{rf:GIKV}, they argued the relation of the homological invariants
with the refined topological vertex \cite{rf:IKV,rf:AK05} and 
the superpolynomial $\overline{\cal P}_{\lambda, \mu} (\bfa, \bq,\bt)$ 
was expressed as a summation over all the partitions (see section 4). 
However, it is totally unclear whether 
the proposed expression actually gives a polynomial in $\bfa$ and 
a specialization $\bfa = \bq^N$ leads the polynomial in $\bq$ and $\bt$.
Moreover, the coefficients of the polynomial after the specialization
are expected to be (non-negative) integers due to the following reason.
%
In \cite{rf:GSV}, it was argued that homological link invariants are
related to a refinement of the BPS state counting in topological open string theory. 
The GIKV conjecture on homological link invariants of the Hopf link 
was based on this proposal. This means for the Hopf link $L$ that 
there is a doubly graded homology theory ${\cal H}_{i,j}^{\mathfrak{sl}(N); \lambda, \mu} (L)$ 
whose graded Poincar\'e polynomial is
\be
\overline{\cal P}_{\mathfrak{sl}(N); \lambda, \mu} (\bq,\bt) = \sum_{i,j \in {\mathbb Z}} 
\bq^i \bt^j \dim {\cal H}_{i,j}^{\mathfrak{sl}(N); \lambda, \mu} (L)~. \label{super1}
\ee
The physical interpretation of ${\cal H}_{i,j}^{\mathfrak{sl}(N); \lambda, \mu} (L)$
as the Hilbert space of BPS states leads to
a prediction on the dependence of the link homologies on the rank $N-1$. 
It has been conjectured \cite{rf:DGR, rf:GW} that there exists a superpolynomial
$\overline{\cal P}_{\lambda, \mu} (\bfa, \bq,\bt)$ 
which is a rational function in three variables such that
\be
\overline{\cal P}_{\mathfrak{sl}(N); \lambda, \mu} (\bq,\bt) =
\overline{\cal P}_{\lambda, \mu} (\bfa ={\bq}^N, \bq,\bt)~. 
\label{super2}
\ee
We see that \eqref{super1} and \eqref{super2} imply that the specialization
 $\bfa = \bq^N$ leads to a polynomial with (non-negative) integer coefficients. 

In this paper, we prove the following. 
\begin{description}
\item{(\romannumeral1)}
The  GIKV's proposal 
$\overline{\cal P}_{\lambda, \mu} ^{\rm GIKV}
(\bfa, \bq,\bt)$ 
for the superpolynomial 
$\overline{\cal P}_{\lambda, \mu} (\bfa, \bq,\bt)$
of the homological invariants of the colored Hopf link
gives really a polynomial in $\bfa$ for arbitrary representations.
\item{(\romannumeral2)}
$\overline{\cal P}_{\lambda, 1^s}^{\rm GIKV} ({\bq}^N, \bq,\bt)$
vanishes for sufficiently small $N\in\bN$,
if one representation is antisymmetric and the other is arbitrary.
\item{(\romannumeral3)}
$\overline{\cal P}_{1^r, 1^s}^{\rm GIKV} ({\bq}^N, \bq,\bt)$ with arbitrary $N\in\bN$
is a polynomial in $\bq$ and $\bt$ with non-negative integer coefficients,
if both representations are antisymmetric.
\end{description}
Furthermore, we perform the summation over the partitions
and show a closed formula 
without assuming the condition $|t|<1$ in \cite{rf:AK09} 
for the superpolynomial of the homological invariants
for antisymmetric representations.
%


The key ingredients of our proof are the Macdonald operators
and a non-standard scalar product for the Hall-Littlewood polynomials.
The Macdonald operators are a set of difference operators
which commute with each other on the space of symmetric functions.
Their  simultaneous eigenfunctions are the Macdonald polynomials
which are two parameter $(q,t)$ deformation of the Schur polynomials.%
\footnote{
We also use another set of parameters $(Q,q,t)$ defined by
$(\bfa; \bq,\bt) = (Q^{-\ha} ;t^{-\ha},-(t/q)^\ha)$.
}{ }
By using a power sum realization for the Macdonald operators
\cite{rf:AKOS,rf:Shi05},%
\footnote{
We use a slightly different definition from that in \cite{rf:Shi05}.
}{ }
we investigate the polynomiality and integrality of 
the superpolynomial of the homological invariants 
obtained by GIKV's proposal.
 

On the other hand, the Hall-Littlewood polynomial is 
the symmetric polynomial obtained from the Macdonald polynomial
by letting $q=0$.
The higher Macdonald operator is realized by using a (non-standard) 
scalar product for which 
the Hall-Littlewood polynomials are 
the basis dual to dominant monomials \cite{rf:DL}.
We show the pairwise orthogonality of  the Hall-Littlewood polynomials
for this scalar product 
and use it for the calculation of the homological invariants.


For non-antisymmetric representations,
the GIKV's proposal for the superpolynomial of the homological invariants 
calls for some improvement 
because it has negative integer coefficients in general.
We find that
this positivity problem may be overcome 
by replacing the Schur function by the Macdonald function 
with an appropriate specialization.




The paper is organized as follows.
In section \ref{sec:HL}, 
we discuss two kinds of scalar product for the Hall-Littlewood polynomials,
which are used in the following sections.
In section \ref{sec:MO}, we give a power sum realization for the Macdonald operators 
which are defined in appendix {\appMacOp}.
Our main results are stated and proved in section \ref{sec:GIKV}. 
We review the Gukov {\it et al}'s conjectures 
on the homological invariants of the colored Hopf link
and prove some of them.
Section \ref{sec:Positivity} and appendix {\appPositivity} are 
devoted to the discussion on the positivity problem and its example, respectively.
Appendix {\appSymm} contains a brief summary of the symmetric functions.  
Proofs of the key theorems are shown in appendices
{\appProofHL} and {\appProofMO}.
In appendix {\appTorus},
we comment on a relation between the torus knot and the Macdonald polynomial.

\vskip12pt
\noindent{\bf Notations.} 
%
The following notations are used through this paper.
%
%
%
Let $\lambda$ be a Young diagram,
i.e. 
a partition $\lambda = (\lambda_1,\lambda_2,\cdots)$,
which is a sequence of non-negative integers such that
$\lambda_{i} \geq \lambda_{i+1}$ and 
$|\lambda| := \sum_i \lambda_i < \infty$.
$\lambda^\vee $ is its conjugate (dual) diagram.
$\ell(\lambda) := \lambda^\vee_1$ is the length. 
%
%
Let $p=(p_1,p_2,\cdots)$ be the power sum symmetric functions in 
\hbox{$x=(x_1,x_2,\cdots)$}
defined as $p_n(x):=\sum_{i\geq 1} x_i^n$.
We treat any symmetric function $f(x)$ in $x$ 
as a function in $p$ unless otherwise stated, 
and sometimes denote it as $f(x(p))$.
$P_\lambda(x;q,t)$, $P_\lambda(x;t)$, 
$s_\lambda(x)$ and $e_\lambda(x)$ are 
the Macdonald, the Hall-Littlewood, the Schur and the elementary symmetric function
in $x$, respectively.
We define the following specialization:
\be
p_n(q^\lambda t^\rho) 
:=
\sum_{i=1}^{\ell(\lambda)} (q^{n\lambda_i}-1)t^{n(\ha-i)}
+ {1 \over t^{n\over 2} - t^{-{n\over 2}} }
=
\sum_{i=1}^N q^{n\lambda_i}t^{n(\ha-i)}
+ {t^{-nN} \over t^{n\over 2} - t^{-{n\over 2}} },
\ee
which is independent of $N$ for any $N\geq \ell(\lambda)$. 
We do not have to assume that $|t|>1$.
Let
$p_n(x,y) := p_n(x) + p_n(y)$;
then, 
\be
p_n(c q^\lambda t^\rho, L t^{-\rho} ) 
=
c^n\sum_{i=1}^{\ell(\lambda)} (q^{n\lambda_i}-1)t^{n(\ha-i)}
+ {c^n - L^n \over t^{n\over 2} - t^{-{n\over 2}} },
\qquad 
c,L \in\bC. 
\ee
%
For $N\in\bZ$ and $N\geq\ell(\lambda )$,
$
p_n(q^\lambda t^{N+\rho}, t^{-\rho} ) 
=
\sum_{i=1}^N q^{n\lambda_i}t^{n(N+\ha-i)}
$
is the power sum in $N$ variables
$\{ q^{\lambda_i} t^{N+\ha-i} \}_{1\leq i \leq N}$.
%
%
We do not have to use any analytic continuation or approximation
because we treat any functions as formal power series. 
For example,
\break\noindent
$ 
\Pi(x,y;q,t)
:=
\Exp{
\sum_{n>0}
{1\over n} 
{1-t^n\over 1-q^n}
p_n(x) p_n(y)
}
$ and 
$
\Dp(x;t) 
:= 
\Exp{-\sum_{n>0}{1-t^{n}\over n}\sum_{i<j}{x_{i}^n\over x_{j}^n} }
$ 
are power series in formal variables $p_n(x)$, $p_n(y)$ and $x_i$
respectively.
%
%
We also let 
$
\gqt{\lambda } 
:=
\prod_{(i,j)\in\lambda }(-1)
q^{\lambda _i-j}t^{-\lambda ^\vee_j+i}
$
and 
$v :=\left({q/t}\right)^{\ha}$.
%
%
The $q$-integer  
${\qint Nt} := {1-t^N\over 1-t}$
and the $q$-binomial coefficient 
${{\qbin Nrt} :=\prod_{i=1}^r{1-t^{N-r+i}\over 1-t^i} }$
are polynomials in $t$ with non-negative integer coefficients for $N,r\in\bN$.

\Section{Scalar product for the Hall-Littlewood polynomial}
\label{sec:HL}%

First, we discuss two kinds of scalar product for symmetric polynomials,
which we will use later.
In this section, 
we treat symmetric functions as symmetric polynomials in finite number of variables 
$x=(x_1,\cdots,x_N)$ by setting $x_i=0$ for $i\geq N+1$. 
For a partition $\lambda=(\lambda _1,\cdots,\lambda_N)$,
we denote
$x^\lambda := x_1^{\lambda _1} x_2^{\lambda _2} \cdots x_N^{\lambda _N}$ and
$m_j := \#\{\ \lambda_i\ |\ \lambda_i=j\ \}$,
i.e.
$\lambda = (0^{m_0}1^{m_1}2^{m_2}\cdots)$ with
$\sum_{i\geq 0} m_i = N$.
The Hall-Littlewood polynomials
$P_\lambda (x;t)$ 
are obtained from the Macdonald polynomials
$P_\lambda (x;q,t)$ defined in appendix {\appSymm} by letting $q=0$:
\ba
P_\lambda (x;t) 
&:=& 
P_\lambda (x;0,t)
\cr
&=&
v_\lambda ^{-1}(t) \sum_{\sigma \in\cS _N}
\sigma \left(x^\lambda \Dp(\bar x;t)^{-1} \right),
\qquad
v_\lambda (t) 
:=
\prod_{j\geq 0} {\qintf {m_j}t}
\label{eq:defHL}%
\cr
\Dp(x;t) 
&:=& 
\Exp
{
-\sum_{n>0}{1-t^{n}\over n}
\sum_{i<j}{x_{i}^n\over x_{j}^n} 
}
=
\prod_{i<j}
\left( 1- {x_i \over x_j }\right)\sum_{n\geq 0}\left({tx_i \over x_j}\right)^n
\ea
with the symmetric group $\cS_N$ on $N$ variables.%
\footnote{
Note that 
$
\Dp(x;t) 
=
\prod_{i < j}
{1-x_i /x_j \over 1-tx_i /x_j }
$
by the analytic continuation.
}{ }
Here
$\bar x := (x_N,x_{N-1},\cdots,x_1)$,
${\qintf Nt}  := {\qint 1t} {\qint 2t} \cdots {\qint Nt}$
and 
${\qint Nt}:={(1-t^n)/(1-t)}$.
The canonical scalar product 
$\langle *,*\rangle_t$ 
is defined from
$\langle *,*\rangle_{q,t}$ 
in appendix {\appSymm} by 
\be
\langle P_\lambda (x;t), P_\mu (x;t)\rangle_t 
:= 
\langle P_\lambda (x;0,t), P_\mu (x;0,t)\rangle_{0,t}
= 
\delta_{\lambda ,\mu} 
\prod_{j\geq 1}\prod_{i=1}^{m_j}{1 \over 1-t^i}.   
\label{eq:scalarHL}%
\ee
We abbreviate it to $\langle P_\lambda , P_\mu\rangle_t$. 
The Cauchy formula \eqref{eq:AppCauchy}
is now
\ba
\sum_\lambda {1\over \langle P_\lambda ,P_\lambda \rangle_t } 
P_\lambda (x;t) P_\lambda (y;t)
=
\Exp{
\sum_{n>0}
{1-t^n\over n}
p_n(x) p_n(y)
}
\label{eq:CauchyHL}%
\ea
which is summed over all partitions $\lambda$.
%


For functions $f$ and $g$ in $x$, 
let us define a second scalar product as \cite{rf:DL}
\be
\langle f(x) , g(x)\rangle_{N;t}'' 
:= 
\oint\prod_{j =1}^N {dx_j \over 2\pi i x_j } 
f(\bar x^{-1}) \Dp(x;t) g(x).
\label{eq:anotherScalarHLdef}%
\ee
Here $x_j$'s  are formal parameters,
and $\exp f$ is defined as the Taylor expansion in $f$. 
For a Laurent series $f(x)$ in $x$,
${\ds\oint{dx\over 2\pi ix} f(x)}$
denotes the constant term in $f$,
i.e.
${\ds\oint{dx\over 2\pi ix} }\sum_{n\in\bZ} f_n x^n = f_0$.
Note that the kernel function $\Dp(x;t)$ is not symmetric in $x_i$. 
Then, 
$
\langle f(x) , g(x)\rangle_{N;t}'' = \langle g(\bar x) , f(\bar x)\rangle_{N;t}''
$.
Hence, the second scalar product is symmetric only for the symmetric functions.
The Hall-Littlewood polynomials $P_\lambda (x;t)$ with $\ell(\lambda )\leq N$
are pairwise orthogonal for the second scalar product and we have
\footnote{
It was shown in \cite{rf:DL} that
the Hall-Littlewood polynomials are the basis dual to dominant monomials.
}{ }
\\{\bf Theorem~2.1.}~~{\it
\be
\langle P_\lambda (x;t) , P_\mu(x;t)\rangle_{N;t}''
=
\langle P_\mu (x;t) , P_\lambda (x;t)\rangle_{N;t}''
=
\delta_{\lambda ,\mu}
v_\lambda (t)^{-1} {\qintf Nt} .
\label{eq:anotherScalarHL}%
\ee
}

A proof is given in appendix {\appProofHL}.
Since 
$P_{1^r} (x;t)$ coincides with the elementary symmetric polynomial $e_r(x)$, 
we obtain
$
\langle P_\lambda (x;t), e_r(x)\rangle_{N;t}''
=
\delta_{\lambda, 1^r} 
{\qbin Nrt},
$
with 
$
{\qbin Nrt}
:=
\prod_{i=1}^r
{ 1-t^{N-r+i} \over 1-t^i }
$.
%
Note that from \eqref{eq:SpecializationE} it follows that
\be
{
\langle e_r, e_r\rangle_{N;t}''
\over
\langle e_r , e_r\rangle_t
}
=
(-1)^r t^{{r\over 2}}
{
e_{N-r}(t^\rho)
\over
e_N(t^\rho)
},
\qquad
e_r(t^\rho)
=
(-1)^r t^{r\over 2}\langle e_r , e_r\rangle_t
,
\label{eq:ScalarScalarHLele}%
\ee
which we will use later.


\Section{Power sum realization for Macdonald operators}
\label{sec:MO}%

In this section, we give a power sum (boson) realization
for the Macdonald operators 
which are defined in appendix {\appMacOp}.
Let $p=(p_1,p_2,\cdots)$ be the power sum symmetric functions 
in infinite number of variables 
$x=(x_1,x_2,\cdots)$ defined as $p_n(x):=\sum_{i\geq 1} x_i^n$.
We treat any symmetric function $f(x)$ in $x$ 
as a function in $p$ unless otherwise stated. 

\subsection{Macdonald operators by the power sum}

%
%
First, we define a commutative family of the difference operators in power sums, 
whose eigenfunctions are the Macdonald functions.
Let $\z=(\z_1,\cdots,\z_r)$.
For $r=0,1,2,\cdots$, let $H$ and $H^r$ be 
$H:=\sum_{r\geq 0}w^r H^r$, $H^0:=1$ and 
\newpage
\ba
H^r &:=&
e_r(t^\rho)
\oint
\prod_{\alpha =1}^r{d\z_\alpha \over 2\pi i\z_\alpha }
\Dp(z;t^{-1}) \vertex^r(\z),
\qquad
\vertex^r(\z)
:=
\vertex^r_+(\z)
\vertex^r_-(\z),
\cr
\vertex^r_+(\z)
&:=&
\Exp
{\sum_{n>0}{1-t^{-n}\over n} 
\sum_{\alpha =1}^r z_\alpha ^n 
p_n },
\cr
\vertex^r_-(\z)
&:=&
\Exp{ \sum_{n>0}{1-t^{-n}\over n}t^n 
\sum_{\alpha =1}^r z_\alpha ^{-n} 
p_n^* }
=
\Exp{ \sum_{n>0}(q^n-1) 
\sum_{\alpha =1}^r z_\alpha ^{-n} 
{\partial\over \partial p_n} }
,
\ea
with
$ e_r(t^\rho) = \prod_{i=1}^r {t^\ha\over t^i-1} $ from \eqref{eq:SpecializationE}
and
$ p_n^* := {1-q^n \over 1-t^n} n {\partial\over \partial p_n } $.
Here $\z_\alpha $ and $p_n$ are formal variables and 
$ \oint{d\z\over 2\pi i\z} f(\z) $ denotes the constant term in $f$.
Note that $H^r$ is written by the second scalar product 
\eqref{eq:anotherScalarHLdef} 
 as
$
H^r = e_r(t^\rho) \langle \vertex^r_+(\z^{-1}),\vertex^r_-(\z)\rangle_{r;t^{-1}}'' 
$.
We also denote them as 
$H(x)$, $H^r(x)$ and $\vertex(\z;x)$
if they act on the power sums $p_n(x)$ in $x$,
but they are independent of the number of variables $x_i$.
%
%

%
%
The Cauchy formula \eqref{eq:CauchyHL} for the Hall-Littlewood function
leads to
\be
\vertex^r(\z) =
\sum_{\lambda \atop \ell(\lambda )\leq r }
{P_\lambda (x(p);t^{-1}) P_\lambda(\z ;t^{-1}) 
\over \langle P_\lambda ,P_\lambda \rangle_{t^{-1}} }
\sum_{\mu \atop \ell(\mu )\leq r }
t^{|\mu |}
{P_\mu (\z^{-1};t^{-1}) P_\mu(x(p^*) ;t^{-1}) 
\over \langle P_\mu ,P_\mu \rangle_{t^{-1}} }.
\label{eq:PhibyHL}%
\ee
Here $P_\lambda (x(p);t)$ is the Hall-Littlewood function in terms of the power sums $p$
and $\ell(\lambda ):=\lambda^\vee_1$ is the length of $\lambda $. 
Thus, $H^r$ is also realized by the Hall-Littlewood function as
\be
H^r =
e_r(t^\rho)
\sum_{\lambda \atop \ell(\lambda )\leq r }
P_\lambda(x(p);t^{-1}) 
{t^{|\lambda |}\langle P_\lambda ,P_\lambda\rangle_{r;t^{-1}}''
\over \langle P_\lambda ,P_\lambda \rangle_{t^{-1}}^2 }
P_\lambda (x(p^*) ;t^{-1}).
\label{eq:HrbyHL}%
\ee
Then, we have 
\cite{rf:Shi05}(\cite{rf:AKOS} for $r=1$)
\footnote{
Our definition is slightly different from \cite{rf:Shi05}.
Ours is not symmetric in $z_\alpha $'s,
and also the integration contour may be different.   
Thus, we give our proof in appendix {\appProofMO}.
}{ }
\\
{\bf Theorem~3.1.}~~{\it
The Macdonald function $P_\lambda(x;q,t)$ is an eigenfunction for $H^r$:
\ba
H P_\lambda(x;q,t) &=& P_\lambda(x;q,t) E_{\lambda},
\cr
H^r P_\lambda(x;q,t) &=& P_\lambda(x;q,t) e_r(q^\lambda t^\rho),
\cr
&&\hskip-100pt
E_{\lambda} 
:=
\Exp{ -\sum_{n>0}{(-w)^n\over n} p_n(q^\lambda t^\rho) }
=
\sum_{r\geq 0} w^r e_r(q^\lambda t^\rho),
\label{eq:EigenH}%
\ea
and therefore, $H^r$ commute with each other $[H^r,H^s]=0$
on the space of symmetric functions. 
}

The proof is given by comparing $H^r$ with the Macdonald operators.

\subsection{Action on the Cauchy kernel}


%
Second, we discuss the properties of $H^r$,
which we will use in subsection \ref{sec:GIKV}.3.
For the Cauchy kernel  $\Pi(x,y;q,t)$ in \eqref{eq:AppCauchy},
we obtain
\be
{\vertex^r(\z;x) \Pi(x,y;q,t) \over \Pi(x,y;q,t) }
=
\Exp
{ \sum_{n>0}{1-t^{-n}\over n}
\left(p_n(x) p_n(z) 
+p_n(ty) p_{-n}(z) 
\right)
},
\ee
with $ p_n(\z) = \sum_{\alpha =1}^r \z_{\alpha }^n $,
($n\in\bZ$).
Since 
\be
\Dp(\z;t^{-1})
\vertex^r(\z ;x) \Pi(x,y;q,t) 
=
\Dp(t\bar\z^{-1};t^{-1})
\vertex^r(t\bar\z^{-1};y) \Pi(x,y;q,t),
\ee
with $\bar\z := (z_r,z_{r-1},\cdots,z_1)$,
the following important duality holds for any variables 
$x=(x_1,x_2,\cdots)$ and $y=(y_1,y_2,\cdots)$,
even though their numbers of components are different:
\be
H^r(x) \Pi(x,y;q,t) 
=
H^r(y) \Pi(x,y;q,t).
\ee

In a way similar to \eqref{eq:PhibyHL} and \eqref{eq:HrbyHL},
we obtain
\be
{\vertex^r(\z;x) \Pi(x,y;q,t) \over \Pi(x,y;q,t) }
=
\sum_{\lambda \atop \ell(\lambda )\leq r }
{P_\lambda (x(p);t^{-1}) P_\lambda(\z ;t^{-1}) 
\over \langle P_\lambda ,P_\lambda \rangle_{t^{-1}} }
\sum_{\mu \atop \ell(\mu )\leq r }
{P_\mu (\z^{-1};t^{-1}) P_\mu(ty(p) ;t^{-1}) 
\over \langle P_\mu ,P_\mu \rangle_{t^{-1}} },
\ee
\be
{H^r(x) \Pi(x,y;q,t) \over \Pi(x,y;q,t) }
=
e_r(t^\rho)
\sum_{\lambda \atop \ell(\lambda )\leq r }
P_\lambda(x(p);t^{-1}) 
{\langle P_\lambda ,P_\lambda\rangle_{r;t^{-1}}''
\over \langle P_\lambda ,P_\lambda \rangle_{t^{-1}}^2 }
P_\lambda (ty(p) ;t^{-1}). 
\label{eq:HrPi}%
\ee


When $x=t^{-\rho}$, we have\footnote{
This proposition was proved in \cite{rf:AK09}
by assuming $|t|<1$.
However, here we do not have to assume it.}
\\
{\bf Proposition~3.2.}~~{\it
\ba
\left.
{H^r(x) \Pi(x,y;q,t)  
\over 
\Pi(x,y;q,t)}
\right\vert_{x=t^{-\rho}}
&=&
\gqt{1^r}
e_r(y,t^{-\rho}),
\label{eq:HPiPi}%
\ea
with
$
\sum_{r\geq 0} (-\w)^r e_r(x,y)
:=
\Exp{-\sum_{n>0}p_n(x,y) \w^n/n }
$
and
$\gqt{1^r}:=\nega{(-1)^r} t^{-{r(r-1)\over 2}}$.
}

\proof 
{}From (\ref{eq:HrPi}) and 
$P_\lambda (t^\rho;t) = \delta_{\lambda ,1^r} e_r(t^\rho)$, 
it follows that
\be
\lhs
=
e_r(t^\rho)
\sum_{s=0}^r 
P_{1^s}(t^{-\rho};t^{-1})
{\langle P_{1^s},P_{1^s}\rangle_{r;t^{-1}}''
\over \langle P_{1^s},P_{1^s}\rangle_{t^{-1}}^2 }
P_{1^s} (ty(p) ;t^{-1}).
\ee
But from (\ref{eq:ScalarScalarHLele}) and \eqref{eq:SpecializationE}, 
we obtain
\be
\lhs
=
\nega{(-1)^r} \prod_{i=1}^r {t^{1-i}}
\sum_{s=0}^r 
e_{r-s}(t^{-\rho})e_s(y).
\ee
Then, 
$\sum_{s=0}^r e_{r-s}(x) e_s(y)=e_r(x,y)$
proves the proposition.
\qed

\subsection{Product of Macdonald operators}


Next we consider the product of the $H^r$'s
which we will use in subsection \ref{sec:GIKV}.4.
For a partition $\lambda $, 
let 
$\ell := \ell(\lambda )$,
$H^\lambda := H^{\lambda _1} H^{\lambda_2} \cdots H^{\lambda_\ell}$,
$\pz^{i } = (\pz^{i }_1,\pz^{i }_2,\cdots,\pz^{i }_{\lambda_i })$
and 
$\tz
:=
(\tz_1,\cdots,\tz_{|\lambda |})
:=
(
\pz^1_1,\cdots,\pz^1_{\lambda_1},
\pz^2_1,\cdots,\pz^2_{\lambda_2},
\cdots,
\pz^\ell_1,\cdots,\pz^\ell_{\lambda_\ell}
)
$;
then,
\be
\Dp(\tz;t) 
=
\Exp{
-\sum_{n>0}{1-t^{n}\over n}
\sum_{i < j  } p_n(\pz^i ) p_{-n}(\pz^j  ) 
}
\prod_{i =1}^\ell \Dp(\pz^i ;t).
\ee
By the OPE relations (as difference operators)
\be
\vertex^s(w)
\vertex^r(z)
=
\Exp{\sum_{n>0}{(1-t^{-n})(q^n-1)\over n} p_{-n}(w) p_n(z)}
\vertex^s_+(w)
\vertex^r_+(z)
\vertex^s_-(w)
\vertex^r_-(z)
\ee
with $z=(z_1,\cdots,z_r)$ and $w=(w_1,\cdots,w_s)$,
we obtain
\ba
&&\hskip-12pt
{
\vertex^{\lambda_\ell}(\pz^{\ell};x)
\cdots
\vertex^{\lambda_1}(\pz^{1};x)
\Pi(x,y;q,t)
\over 
\Pi(x,y;q,t)
}
\cr
&=&
\Exp{
\sum_{n>0}{1-t^{-n}\over n}
\left\{
p_n(x) p_n(\tz)
+
p_n(ty) p_{-n}(\tz)
+ 
(q^n-1)\sum_{i < j  } p_n(\pz^i ) p_{-n}(\pz^j  ) 
\right\} 
}
\cr
&=&
\Exp{
\sum_{n>0}{1-t^{-n}\over n}
\left\{
p_n(x) p_n(\tz)
+
p_n(ty) p_{-n}(\tz)
+ 
\sum_{i < j  } p_n(q \pz^i ) p_{-n}(\pz^j  ) 
\right\} 
}
\cr
&&\hskip250pt
\times
{ \Dp(\tz;t^{-1}) \over \prod_{i =1}^\ell \Dp(\pz^i;t^{-1}) }
.~~~
\ea


When $(x,y)=(t^\rho,ct^{-\rho})$ with $c\in\bC$, 
since
$p_n(t^\rho) = (t^{{n\over 2}}-t^{-{n\over 2}})^{-1}$
it follows that
\ba
&&
\prod_{i =1}^\ell \Dp(\pz^i;t^{-1})
\times
\left.
{
\vertex^{\lambda_\ell}(\pz^{\ell};x)
\cdots
\vertex^{\lambda_1}(\pz^{1};x)
\Pi(x,ct^{-\rho};q,t)
\over 
\Pi(x,ct^{-\rho};q,t)
}
\right\vert_{x=t^\rho}
\cr
&&\hskip24pt
=
\Dp(\tz;t^{-1})
\Exp{
\sum_{n>0}{1-t^{-n}\over n}
\sum_{i < j  } p_n(q \pz^i ) p_{-n}(\pz^j  ) 
}
\prod_{\alpha =1}^{|\lambda |}
{1-t^\ha c/\tz_\alpha \over 1-t^{-\ha} \tz_\alpha },
\ea
which is a polynomial of degree $|\lambda |$ in $c$.
For abbreviation, we write ${1/( 1-z )}$ instead of  $\sum_{n\geq 0}z^n$. 
Therefore, by taking the constant term in $\tz$, we have
\\
{\bf Proposition~3.3.}~~{\it 
\ba
&&
\left.
{
H^{\lambda}(x)
\Pi(x,ct^{-\rho};q,t)
\over 
e_{\lambda}(t^\rho)
\Pi(x,ct^{-\rho};q,t)
}
\right\vert_{x=t^\rho}
\cr
&&\hskip24pt
=
\oint\prod_{\alpha =1}^{|\lambda |} 
{d\tz_\alpha \over 2\pi i \tz_\alpha } 
{1-c/\tz_{\alpha } \over 1-\tz_{\alpha } }
\prod_{\alpha <\beta }^{|\lambda |} 
{1- \tz_\alpha/\tz_\beta \over 1-\tz_\alpha/t \tz_\beta }
\prod_{i < j }^\ell \prod_{\alpha =1}^{\lambda_i } \prod_{\beta =1}^{\lambda _ j }
{1-q\pz^i_\alpha/\pz^j _\beta \over 1-q\pz^i_\alpha/t \pz^j _\beta }
~~~~~
\label{eq:HHPi}%
\ea
is a polynomial of degree $|\lambda |$ in $c$
and a polynomial in $q$ and $1/t$ with integer coefficients
and vanishes when $c=1$.
}

\proof
The \rhs\ 
of \eqref{eq:HHPi} reduces to
\ba
\oint\prod_{\alpha =1}^{|\lambda |} 
{d\tz_\alpha \over 2\pi i \tz_\alpha } 
\left(1-{c\over \tz_{\alpha }}\right)
\sum_{n =0}^\alpha (\tz_\alpha )^n  
&\times&
\prod_{\alpha <\beta }^{|\lambda |} 
\left(1- {\tz_\alpha\over \tz_\beta }\right) 
\sum_{n =0}^\alpha \left({\tz_\alpha\over t \tz_\beta }\right)^n
\cr
&\times&
\prod_{i < j }^\ell \prod_{\alpha =1}^{\lambda_i } \prod_{\beta =1}^{\lambda _ j }
\left(1-{q\pz^i_\alpha \over \pz^j _\beta }\right) 
\sum_{n =0}^{|\pz^i_\alpha |}
\left({q\pz^i_\alpha \over t \pz^j _\beta }\right)^n
,
\ea
with $|\pz^i_\alpha |:=\alpha + \sum_{k=1}^{i-1}\lambda_k $,
so this is a polynomial in $q$ and $1/t$ with integer coefficients.
When $c=1$, 
since $(1-1/\tz_\alpha )\sum_{n\geq 0} (\tz_\alpha )^n = 1/\tz_\alpha $,
the integrand 
in \eqref{eq:HHPi} 
has no constant term in the last variable $z_{|\lambda |}$;
thus,
\eqref{eq:HHPi} vanishes.
\qed

Let
\be
\widetilde 
Z^{\lambda } 
\prod_{j =1}^{\lambda_1} (1-c/t^{j -1})
:=
\left.
{
H^{\lambda }(x)
\Pi(x,ct^{-\rho};q,t)
\over 
e_{\lambda}(t^\rho)
\Pi(x,ct^{-\rho};q,t)
}
\right\vert_{x=t^\rho}
.
\ee
Then, for example,
\ba
\widetilde
Z^{(3)}
&=&
1,
\cr
\widetilde
Z^{(2,1)}
&=&
(1-cq)-c(1-q)/t^2,
\\
\widetilde
Z^{(1,1,1)}
&=&
(1-cq)
\{
(1-cq^2) 
- c(1-q)(1+2q)/t
- c(1-q)^2/t^2
\}
+ c^2(1-q)^2/t^3
.~~~~ \nonumber
\ea
Direct calculation of 
the functions $\widetilde Z^\lambda $ for $|\lambda |\leq 7$
suggests the following conjecture.
\\{\bf Conjecture~3.4.}~~{\it 
If 
$c=t^N$, $N\in\bZ$ with $0\leq N < \lambda _1$,
then $(\ref{eq:HHPi})$ vanishes.
}

\Section{Homological link invariants from GIKV conjecture }
\label{sec:GIKV}%

In this section, we introduce the 
conjecture on homological link invariants of the colored Hopf link
by Gukov, Iqbal, Koz\c{c}az and Vafa \cite{rf:GIKV}
and prove some of them.%
\footnote{
Our 
$(Q,\GIKVQ; q, t;\bq,\bt) $
equals 
$ (Q\sqrt{q_2/q_1},Q; q_1^{-1}, q_2^{-1};q,t) $ of \cite{rf:GIKV}.
\hfill\break
In \cite{rf:GIKV},
$
{\HopfLink \lambda \mu \bfa\bq\bt}
=
(-\bfa)^{|\lambda |+|\mu|}
\bt^{2|\lambda ||\mu|}
G_{\lambda ,\mu}(-\GIKVQ;\bq,\bt) 
$
and 
$
G_{\lambda ,\mu}(-\GIKVQ;\bq,\bt) 
=
Z^{{\rm inst}}_{\lambda ,\mu}(Q;q,t) 
$
,
where
$(\bfa; \bq,\bt) := (Q^{-\ha} ;t^{-\ha},-(t/q)^\ha)$.
Note that
$\bfa = \bq^N$ and $\bt=-1$ are
equivalent to
$Q=t^N$ and $q=t$, respectively.
}{ }

\subsection{GIKV conjecture}


First, we recall the GIKV conjecture and present our main theorem.
Following \cite{rf:GIKV} let us consider
\be
Z_{\lambda ,\mu}(\Q ;q,t) 
:= 
\sum_\eta 
s_{\lambda}(q^{\eta}t^{\rho})
s_{\mu}(q^{\eta}t^{\rho})
\prod_{(i,j)\in \eta} 
{ \vQ
\over 
(1-q^{ \eta_i-j+1}t^{ \eta^\vee_j-i  })
(1-q^{-\eta_i+j  }t^{-\eta^\vee_j+i-1})
}.
\label{eq:GIKVZ}
\ee
Here $s_{\lambda}(q^{\eta}t^{\rho})$ is the Schur function 
in the power sum $p_n(q^{\eta}t^{\rho})$.
{}From \eqref{eq:PPSpecialization} and \eqref{eq:Pdual}, 
it is written by the Macdonald functions 
$P_\lambda (x;q,t)$ defined in appendix {\appSymm} 
as follows:
\ba
Z_{\lambda ,\mu}(\Q ;q,t) 
&=& 
\sum_\eta (-\viQ)^{|\eta|}
P_{\eta^\vee}(q^{\rho};t,q)
P_{\eta}(t^{\rho};q,t)
s_{\lambda}(q^{\eta}t^{\rho})
s_{\mu}(q^{\eta}t^{\rho})
\cr
&=&
\sum_\eta {\Qv^{|\eta|}\over \langle P_\eta , P_\eta \rangle_{q,t} }
P_{\eta}(t^{-\rho};q,t)
P_{\eta}(t^{\rho};q,t)
s_{\lambda}(q^{\eta}t^{\rho})
s_{\mu}(q^{\eta}t^{\rho}).
\ea
Note that 
$
Z_{\bullet,\bullet}(\Q ;q,t) 
=
\Pi(\vQ t^{\rho},t^{-\rho};q,t).
$
Then, 
$
Z^{{\rm inst}}_{\lambda ,\mu}(\Q ;q,t) 
:=
Z_{\lambda ,\mu}(\Q ;q,t) 
/
Z_{\bullet,\bullet}(\Q ;q,t) 
$
is a power series in $\Q $ and a meromorphic function in $q$ and $t$.
But $Z^{{\rm inst}}_{\lambda ,\mu}(\Q ;q,t)$ is expected to give 
the superpolynomial of the homological invariants 
and  Gukov, Iqbal, Koz\c{c}az and Vafa claimed the following results.
\\
{\bf Conjecture~4.1}\cite{rf:GIKV}.~{\it

\noindent {\rm (1)}
$Z^{{\rm inst}}_{\lambda ,\mu}(\Q ;q,t)$
is a finite polynomial in $\Q $.

\noindent {\rm (2)}
$Z^{{\rm inst}}_{\lambda ,\mu}(\vitN ;q,t)$
vanishes for sufficiently small $N\in\bZ_{\geq 0}$.

\noindent {\rm (3)}
$Z^{{\rm inst}}_{\lambda ,\mu}(\vitN ;q,t)\nega{(-1)^{|\lambda |+|\mu|}}$
$(N\in\bZ_{\geq 0})$ is a finite 
polynomial in $q$ and $t$
with integer coefficients.

\noindent {\rm (4)}
$Z^{{\rm inst}}_{\lambda ,\mu}(\vitN ;q,t)\nega{(-1)^{|\lambda |+|\mu|}}$
for sufficiently large $N$
coincides with the $sl(N)$  homological invariants of the Hopf link
colored by the representations $\lambda $ and $\mu$
up to an overall factor.
}

Note that if the conjecture 4.1.(1) is established, one can easily calculate  
$Z^{{\rm inst}}_{\lambda ,\mu}(\vitN ;q,t)$
for any given  $(\lambda ,\mu; N)$ 
and check (2) and (3).
In the following subsections, 
we will show part of the GIKV conjecture by proving a series of propositions.%
\footnote{
In \cite{rf:AK09},
we have shown 
(\romannumeral2)--(\romannumeral4) by assuming the condition $|t|<1$.
}{ }
%
\\
{\bf Theorem~4.2.}~~{\it

\noindent
$(\romannumeral1)$ Conjecture
$(1)$ holds for arbitrary $(\lambda ,\mu)$
$($proposition $4.10)$.

\noindent
$(\romannumeral2)$ Conjecture
$(2)$ holds for arbitrary $(\lambda ,\mu)$ with $|\lambda |+|\mu|\leq 7$ or $\mu =1^s$
$($propositions $4.10$ and $4.8)$.

\noindent
$(\romannumeral3)$ Conjecture
$(3)$ holds for  $(\lambda ,\mu) = (\lambda ,1^s)$ 
with $|\lambda |\leq 7$ or $\lambda =1^r$
$($propositions $4.8$ and $4.6)$.

\noindent
$(\romannumeral4)$ Conjecture
$(4)$ holds for  $(\lambda ,\mu) = (1^r ,1)$
$($proposition $4.6)$.
}
\newcommand{\propMainTheorem}{4.2}

\subsection{$q=t$ case}


When $q=t$,  we have
\\
{\bf Proposition~4.3.}~~{\it 
\be
{
Z_{\lambda,\mu}(Q ;q,q) 
\over 
Z_{\bullet,\bullet}(Q ;q,q) 
}
=
\gqq{\lambda}\gqq{\mu}
s_\lambda(Q q^\rho,q^{-\rho})
s_\mu(Q q^{\lambda} q^\rho, q^{-\rho}),
\ee
which is a polynomial of degree $|\lambda|+|\mu|$ in $Q $. 
Here
$\gqq{\lambda}:=\prod_{(i,j)\in\lambda}(-1) q^{\lambda_i-\lambda^\vee_j+i-j}$.
}

\proof
First, we have the following nontrivial identity:
\be
s_\lambda(q^\rho)s_\mu(q^{\lambda+\rho})  
=
\gqq{\lambda} \gqq{\mu}
\sum_\nu
s_{\lambda/\nu}(q^{-\rho})s_{\mu/\nu}(q^{-\rho})  
=
s_\mu(q^\rho)s_\lambda(q^{\mu+\rho}),  
\label{eq:cyclic}
\ee
which is proved by 
$s_{\lambda/\mu}(q^\rho)= 
s_{\lambda^\vee/\mu^\vee}(-q^{-\rho})$
and the cyclic symmetry of the topological vertex \cite{AKMV}.
Since the above identity and 
$
s_{\lambda}(q^\rho)
= 
\gqq{\lambda}
s_{\lambda^\vee}(-q^\rho),
$
it follows that
\be
Z_{\lambda,\mu}(\Q ;q,q) 
= 
\gqq{\mu}
\sum_{\eta,\nu} Q^{|\eta|}
s_{\mu/\nu}(q^{-\rho})
s_{\eta/\nu}(q^{-\rho})
s_\lambda(q^\rho)
s_\eta(q^{\lambda+\rho}).
\ee
The Cauchy formula \eqref{eq:AppCauchy}
and the adding formula \eqref{eq:AppaddSkewMacdonald}
yield
\be
Z_{\lambda,\mu}(\Q ;q,q) 
= 
\gqq{\mu}
s_\mu(Qq^{\lambda+\rho},q^{-\rho})
\Pi(Qq^{\lambda+\rho},q^{-\rho};q,q)
s_{\lambda}(q^\rho).
\ee
But ((2.12) and (5.20) of \cite{rf:AK08})
\be
{
\Pi(Q q^{\lambda+\rho},q^{-\rho};q,q)
\over 
\Pi(Q q^{\rho},q^{-\rho};q,q)
}
=
\gqq{\lambda}
{
s_\lambda (Q q^{\rho},q^{-\rho})
\over 
s_\lambda (q^{\rho})},
\ee
which completes the proof.
\qed


Note that for $N\in\bZ$ and $N\geq\ell(\lambda )$,
$
p_n(q^{\lambda+N+\rho}, q^{-\rho} ) 
=
\sum_{i=1}^N q^{n(\lambda_i+i-\ha)}.
$
Thus,
$s_\mu (q^{\lambda+N+\rho}, q^{-\rho} )$
is the Schur polynomial in $N$ variables
$\{ q^{\lambda_i+i-\ha} \}_{1\leq i \leq N}$,
which is a polynomial in $q$ with non-negative integer coefficients
and vanishes for $\ell(\lambda )\leq N < \ell(\mu)$.
Therefore, when $Q =q^N$, we have 
\\
{\bf Proposition~4.4.}~~{\it
If $N\in\bZ_{\geq 0}$,
\ba
(-1)^{|\lambda|+|\mu|}
{
Z_{\lambda,\mu}(q^N ;q,q) 
\over 
Z_{\bullet,\bullet}(q^N ;q,q) 
}
=
(-1)^{|\lambda|+|\mu|}
\gqq{\lambda} \gqq{\mu}
s_\mu\left(\{ q^{\lambda_i+i-\ha} \}_{1\leq i \leq N}\right)
s_\lambda\left(\{ q^{i-\ha} \}_{1\leq i \leq N}\right),
\ea
which is a polynomial in $q$ with non-negative integer coefficients
and vanishes for $0\leq N < \max \{\ell(\lambda),\ell(\mu)\}$.
This coincides with the colored Hopf link invariant by 
}\cite{rf:ML}{\it
up to the over all factor.%
\footnote{
This was shown in \cite{rf:GIKV} for $N\rightarrow\infty$.
}{ }
}

\proof
If $0\leq N< \ell(\lambda)$, $s_\lambda (q^{N+\rho},q^{-\rho})=0$.
If $\ell(\lambda)\leq N< \ell(\mu)$, $s_\mu(q^{\lambda+N+\rho}, q^{-\rho})=0$.
By the symmetry 
$Z_{\lambda,\mu}(Q;q,q)=Z_{\mu,\lambda}(Q;q,q)$,
we have
$Z_{\lambda,\mu}(Q;q,q)=0$
for 
$0\leq N < \max \{\ell(\lambda),\ell(\mu)\}$.
On the other hand, if $N \geq \max \{\ell(\lambda),\ell(\mu)\}$,
$s_\mu(q^{\lambda+N+\rho}, q^{-\rho})$
is the Schur polynomial in $N$ variables
$\{ q^{\lambda_i+i-\ha} \}_{1\leq i \leq N}$,
which is a polynomial in $q$ with non-negative coefficients. 
\qed

When $q\neq t$, 
it is difficult to calculate $Z_{\lambda,\mu}(Q;q,t)$ explicitly 
for lack of the cyclic symmetry of the refined topological vertex \cite{rf:GIKV}.

\subsection{$(\lambda,\mu)=(1^r,1^s)$ case}


For $\lambda = 1^r$ and $\mu = 1^s$,  we have
\\
{\bf Proposition~4.5.}~~{\it 
\be
{
Z_{1^r,1^s}(\Q ;q,t) 
\over 
Z_{\bullet,\bullet}(\Q ;q,t) 
}
=
\gqt{1^r}\gqt{1^s} 
e_r(\vQ t^\rho,t^{-\rho})
e_s(\vQ q^{(1^r)} t^\rho, t^{-\rho}),
\label{eq:AntiAnti}%
\ee
which is a polynomial of degree $r+s$ in $\Q $. 
Here
$\gqt{1^r}:=\nega{(-1)^r} t^{-{r(r-1)\over 2}}$.
}

\proof
Since
$s_{1^r}(x)= P_{1^r}(x;q,t)= e_r(x)$,
it follows that
\be
Z_{1^r,1^s}(\Q ;q,t) 
= 
\sum_\eta {\Qv^{|\eta|}\over \langle P_\eta , P_\eta \rangle_{q,t} }
P_{\eta}(t^{-\rho};q,t)
P_{\eta}(t^{\rho};q,t)
P_{1^r}(q^{\eta}t^{\rho};q,t)
e_{s}(q^{\eta}t^{\rho}).
\ee
The symmetry 
$
P_\lambda(t^\rho;q,t)P_\mu(q^\lambda t^\rho;q,t)
=
P_\mu(t^\rho;q,t)P_\lambda(q^\mu t^\rho;q,t)
$
in \eqref{eq:Wsymm} leads to 
\be
Z_{1^r,1^s}(\Q ;q,t) 
= 
P_{1^r}(t^{\rho};q,t)
\sum_\eta 
{\Qv^{|\eta|}\over \langle P_\eta , P_\eta \rangle_{q,t} }
P_{\eta}(t^{-\rho};q,t)
P_{\eta}(q^{(1^r)}t^{\rho};q,t)
e_{s}(q^{\eta}t^{\rho}).
\ee
By \eqref{eq:EigenH}, 
we can replace 
$e_{s}(q^{\eta}t^{\rho})$
by $H^s$ as follows:
\be
Z_{1^r,1^s}(\Q ;q,t) 
= 
P_{1^r}(t^{\rho};q,t)
\sum_\eta {1\over \langle P_\eta , P_\eta \rangle_{q,t} }
P_{\eta}(\vQ q^{(1^r)}t^{\rho};q,t)
H^s P_\eta(x;q,t)
\vert_{x=t^{-\rho}}.
\ee
The Cauchy formula \eqref{eq:AppCauchy} and \eqref{eq:HPiPi}
yield
\ba
Z_{1^r,1^s}(\Q ;q,t) 
&=& 
P_{1^r}(t^{\rho};q,t)
H^s
\Pi(x,\vQ q^{(1^r)}t^{\rho};q,t)
\vert_{x=t^{-\rho}}
\cr
&=& 
\gqt{1^s}
P_{1^r}(t^{\rho};q,t)
\Pi(t^{-\rho},\vQ q^{(1^r)}t^{\rho};q,t)
e_s(\vQ q^{(1^r)} t^\rho, t^{-\rho}).
\ea
But ((2.12) and (5.20) of \cite{rf:AK08})
\be
{
\Pi(\vQ q^{(1^r)}t^{\rho},t^{-\rho};q,t)
\over 
\Pi(\vQ t^{\rho},t^{-\rho};q,t)
}
=
\gqt{1^r}
{
P_{1^r} (\vQ t^{\rho},t^{-\rho};q,t)
\over 
P_{1^r} (t^{\rho};q,t)
},
\ee
which completes the proof.
\qed


Note that for $N\in\bZ$ and $N\geq\ell(\lambda )$,
$
p_n(q^\lambda t^{N+\rho}, t^{-\rho} ) 
=
\sum_{i=1}^N q^{n\lambda_i}t^{n(N+\ha-i)}
$.
Thus,
$e_r (q^\lambda t^{N+\rho}, t^{-\rho} )$
is the elementary symmetric polynomial in $N$ variables
$\{ q^{\lambda_i} t^{N+\ha-i} \}_{1\leq i \leq N}$,
which is a polynomial in $q$ and $t$ with non-negative integer coefficients
and vanishes for $\ell(\lambda )\leq N < r$.
Therefore, when $\vQ =t^N$, we have 
\\
{\bf Proposition~4.6.}~~{\it
If $N\in\bZ_{\geq 0}$,
\be
\nega{(-1)^{r+s}}
{
Z_{1^r,1^s}(\vitN ;q,t) 
\over 
Z_{\bullet,\bullet}(\vitN ;q,t) 
}
=
\nega{(-1)^{r+s}}
\gqt{1^s}
\gqt{1^r}
e_r (\{ t^{i-\ha} \}_{1\leq i \leq N})
e_s(\{ q^{(1^r)_i} t^{N+\ha-i} \}_{1\leq i \leq N}),
\label{eq:AntiAntiN}
\ee
which is a polynomial in $q$ and $t$ with non-negative integer coefficients
and vanishes for $0\leq N < \max \{r,s\}$ 
$($theorem $\propMainTheorem$~$(\romannumeral3)$, $\lambda=1^r$ case$)$.
}


When $s=1$, up to an overall factor,
\eqref{eq:AntiAntiN} agrees to the recent result by Yonezawa \cite{rf:Yone} 
who used the method of matrix  factorization
developed by Khovanov and Rozansky \cite{rf:KR}
(theorem $\propMainTheorem$~$(\romannumeral4)$).


In \cite{rf:Taki08},
based on an assumption on the computation of the refined topological vertex, 
Taki proposed that the homological invariants for the colored Hopf link are written as  
\footnote{
Our $(q,t) = (q^{-1},t^{-1}) $ 
of \cite{rf:Taki08}. 
}{ }
\ba
\bar {\cal P}^{{\rm Taki}}_{\mu,\lambda}(\viQ ;q^{-1},t^{-1})  
&:=& 
\bar c_{\lambda,\mu}
s_{\lambda}(t^{\rho})
s_{\mu}(q^{\lambda}t^{\rho},\ \Qv^{-1}t^{-\rho})
N_{\bullet,\lambda}(\Qv^{-1} q/t;q,t) 
\gqt{\lambda}\pos{(-1)^{|\lambda|}}
\cr
&=& 
\bar c_{\lambda,\mu}
s_{\lambda}(t^{\rho})
s_{\mu}(\vQ q^{\lambda}t^{\rho},\ t^{-\rho})
N_{\lambda,\bullet}(\vQ ;q,t),
\cr
\bar c_{\lambda,\mu} 
&:=&
(-1)^{|\lambda|}\Qv^{{|\lambda|+|\mu|\over 2}} v^{-2|\lambda||\mu|}.
\ea
with the denominator factor $N_{\lambda,\mu}$ of Nekrasov's partition function:
\be
N_{\lambda,\mu}(\vQ ;q,t)
:=
\prod_{(i,j)\in\lambda} 
\left( 1 - \vQ\, q^{\lambda_i-j} t^{\mu^\vee_j-i+1} \right)
\prod_{(i,j)\in\mu } 
\left( 1 - \vQ\, q^{-\mu_i+j-1} t^{-\lambda^\vee_j+i  } \right).
\ee
For $(\lambda,\mu) = (1^r,1^s)$,  
\be
\bar {\cal P}^{{\rm Taki}}_{1^s,1^r}(\viQ;q^{-1},t^{-1}) 
= 
\bar c_{1^r,1^s}
s_{1^r}(t^{\rho})
e_s(\vQ q^{(1^r)}t^{\rho},\ t^{-\rho})
N_{1^r,\bullet}(\vQ ;q,t),
\ee
which coincides with (\ref{eq:AntiAnti}) 
up to the over all factor 
$\nega{(-1)^{r+s+1}}\pos{-} \Qv^{(r+s)/2}v^{-2rs}$.

\subsection{$(\lambda,\mu)=(\lambda,1^s)$ case}


Since the set of the Macdonald functions $\{P_\lambda (x;q,t)\}_{|\lambda |=d}$
is a basis of the homogeneous symmetric functions of degree $d$, 
we can write the Schur function $s_{\lambda}(x)$
by the Macdonald functions as
\be
s_{\lambda}(x) = 
\sum_{\mu \leq \lambda }
U_{\lambda,\mu} P_\mu (x;q,t),
\qquad
U_{\lambda ,\mu}
:=
\sum_{\nu (\lambda \geq \nu \geq\mu )} 
u_{\lambda , \nu }(q,q) (u^{-1})_{\nu, \mu}(q,t),
\label{eq:sP}
\ee
where $u_{\lambda , \mu}(q,t)$ is defined by \eqref{eq:MacDef}.
Note that $U_{\lambda ,\mu}$ is a rational function in $q$ and $t$.
Then, we have
\\
{\bf Proposition~4.7.}~~{\it 
\be
{
Z_{\lambda ,1^s}(\Q ;q,t) 
\over 
Z_{\bullet,\bullet}(\Q ;q,t) 
}
=
\gqt{1^s}
\sum_{\mu \leq \lambda }
U_{\lambda,\mu}
\gqt{\mu} 
P_\mu (\vQ t^{\rho},t^{-\rho};q,t)
e_s(\vQ q^\mu  t^\rho, t^{-\rho}),
\ee
which is a polynomial of degree $|\lambda | + s$ in $\Q $. 
Here
$
\gqt{\lambda}
:=
\prod_{(i,j)\in\lambda }\nega{(-1)} q^{\lambda_i-j}t^{-\lambda^\vee_j+i}
$.
}

\proof
First, we have
\be
Z_{\lambda ,1^s}(\Q ;q,t) 
=
\sum_{\mu \leq \lambda }
U_{\lambda,\mu}
\widetilde Z_{\mu ,1^s}(\Q ;q,t), 
\ee
where
\be
\widetilde Z_{\mu ,1^s}(\Q ;q,t) 
:=
\sum_\eta 
{\Qv^{|\eta|}\over \langle P_\eta , P_\eta \rangle_{q,t} }
P_{\eta}(t^{-\rho};q,t)
P_{\eta}(t^{\rho};q,t)
P_{\mu}(q^{\eta}t^{\rho};q,t)
s_{1^s}(q^{\eta}t^{\rho}).
\ee
Then, \eqref{eq:Wsymm} 
leads to 
\be
\widetilde Z_{\mu ,1^s}(\Q ;q,t) 
=
\sum_\eta 
{\Qv^{|\eta|}\over \langle P_\eta , P_\eta \rangle_{q,t} }
P_{\eta}(t^{-\rho};q,t)
P_{\mu}(t^{\rho};q,t)
P_{\eta}(q^\mu t^{\rho};q,t)
s_{1^s}(q^{\eta}t^{\rho}).
\ee
We can proceed in the same way as the last subsection.
Namely, 
$s_{1^s}(x) = e_s(x)$,
(\ref{eq:EigenH}) and 
(\ref{eq:HPiPi}) yield
\ba
\widetilde Z_{\mu ,1^s}(\Q ;q,t) 
&=&
P_{\mu}(t^{\rho};q,t)
\sum_\eta {\Qv^{|\eta|}\over \langle P_\eta , P_\eta \rangle_{q,t} }
P_{\eta}(q^\mu t^{\rho};q,t)
H^s P_{\eta}(x;q,t)
\vert_{x=t^{-\rho}}.
\cr
&=&
P_{\mu}(t^{\rho};q,t)
H^s
\Pi(x,\vQ q^{\mu }t^{\rho};q,t)
\vert_{x=t^{-\rho}}
\cr
&=&
P_{\mu}(t^{\rho};q,t)
\gqt{1^s}
e_s(\vQ q^\mu t^\rho,t^{-\rho})
\Pi(t^{-\rho},\vQ q^{\mu }t^{\rho};q,t).
\ea
But from ((2.12) and (5.20) of \cite{rf:AK08})
\be
{
\Pi(\vQ q^{\mu}t^{\rho},t^{-\rho};q,t)
\over 
\Pi(\vQ t^{\rho},t^{-\rho};q,t)
}
=
\gqt{\mu}
{
P_\mu (\vQ t^{\rho},t^{-\rho};q,t)
\over
P_\mu (t^{\rho};q,t)
},
\ee
we conclude that
\be
{
\widetilde Z_{\mu ,1^s}(\Q ;q,t) 
\over 
Z_{\bullet,\bullet}(\Q ;q,t) 
}
=
\gqt{1^s}\gqt{\mu}
P_\mu (\vQ t^{\rho},t^{-\rho};q,t)
e_s(\vQ q^\mu  t^\rho, t^{-\rho}).
\ee
\qed


Note that for $N\in\bZ$ and $N\geq\ell(\lambda )$,
$P_\mu (q^\lambda t^{N+\rho}, t^{-\rho} )$
is the Macdonald polynomial in $N$ variables
$\{ q^{\lambda_i} t^{N+\ha-i} \}_{1\leq i \leq N}$
and vanishes for $\ell(\lambda )\leq N <\ell(\mu)$.
Therefore, when $\vQ =t^N$, we have 
\\
{\bf Proposition~4.8.}~~{\it
If $N\in\bZ_{\geq 0}$,
\ba
{
Z_{\lambda ,1^s}(\vitN ;q,t) 
\over 
Z_{\bullet,\bullet}(\vitN ;q,t) 
}
&=&
\gqt{1^s}
\sum_{\mu \leq \lambda }
U_{\lambda ,\mu} \gqt{\mu}
P_\mu (\{ t^{i-\ha} \}_{1\leq i \leq N};q,t)
e_s(\{ q^{\mu_i} t^{N+\ha-i} \}_{1\leq i \leq N})
\cr
&=&
\gqt{1^s}
e_s(\{ t^{i-\ha} \}_{1\leq i \leq N})
\sum_{\mu \leq \lambda }
U_{\lambda ,\mu} \gqt{\mu}
P_\mu(\{ q^{(1^s)_i} t^{N+\ha-i} \}_{1\leq i \leq N};q,t),
\label{eq:LaAntiQtN}%
\ea
which vanishes for $0\leq N < \max \{\ell(\lambda ),s\}$ 
$($theorem $\propMainTheorem$~$(\romannumeral2)$, $\mu=1^s$ case$)$.
}


Finally, we should make a remark on the fact that
the transition function $U_{\lambda ,\mu}(q,t)$ in \eqref{eq:sP}
is a rational function in $q$ and $t$. 
It is not obvious that the formulas in the above propositions 
are in fact polynomials in $q$ and $t$. 
However, we have checked the following conjecture up to $d=7$
by direct calculation.
\\
{\bf Conjecture~4.9.}~~{\it
For $|\lambda|=d$,
$
\sum_{\mu ( \lambda \geq \mu \geq \nu ) }
U_{\lambda ,\mu} \gqt{\mu}
(U^{-1})_{\mu, \nu}
$
is a polynomial with integer coefficients 
of degree $d(d-1)/2$ in $q$ and 
of degree $d(d-1)/2$ in $t^{-1}$.
}

On the assumption that the above conjecture is true,
if $N\in\bZ_{\geq 0}$ then
\ba
{
Z_{\lambda ,1^s}(\vitN ;q,t) 
\over
Z_{\bullet,\bullet}(\vitN ;q,t) 
}
&=&
\gqt{1^s}
e_s(t^{N+\rho},t^{-\rho})
\sum_{\mu,\nu\atop |\mu|=|\nu|=|\lambda |} 
U_{\lambda ,\mu} \gqt{\mu}
(U^{-1})_{\mu, \nu}
s_\nu(q^{(1^s)} t^{N+\rho}, t^{-\rho})
\ea
is a polynomial in $q$ and $t$ with integer coefficients
because $s_\lambda (x)$ is a function in $x$ with non-negative integer coefficients
(theorem $\propMainTheorem$~$(\romannumeral3)$, $|\lambda|\leq 7$ case).

\subsection{General $(\lambda,\mu)$ case}


From \eqref{eq:eleUT} and \eqref{eq:MacDef}, we have
\be
s_{\lambda}(x) = 
\sum_{\mu \leq \lambda }
V_{\lambda,\mu} e_{\mu^\vee} (x),
\qquad
V_{\lambda ,\mu}
:=
\sum_{\nu (\lambda \geq \nu \geq\mu )} 
u_{\lambda , \nu } (a^{-1})_{\nu, \mu},
\label{eq:se}
\ee
where $u_{\lambda , \mu}:= u_{\lambda , \mu}(q,q)$ 
and $a_{\lambda , \mu}$ is defined in \eqref{eq:eleUT}.
Note that $|\mu^\vee|=|\lambda |$ in the above equation.
More precisely, we have the Jacobi-Trudy formula 
$
s_{\lambda}(x) 
= 
\det\left(e_{\lambda^\vee_i - i + j}(x)\right)_{1\leq i,j \leq \lambda_1}
$
with $e_{-r}=0$ for $r>0$.
Then, we have
\\
{\bf Proposition~4.10.}~~{\it 
\be
{
Z_{\lambda ,\mu}(\Q ;q,t) 
\over 
Z_{\bullet,\bullet}(\Q ;q,t) 
}
=
\sum_{\nu\leq\lambda \atop \sigma\leq\mu }
V_{\lambda ,\nu} V_{\mu,\sigma}
\left.
{
H^{\nu^\vee} (x)H^{\sigma^\vee} (x)
\Pi(x,\vQ t^{-\rho};q,t)
\over 
\Pi(x,\vQ t^{-\rho};q,t)
}
\right\vert_{x=t^\rho}, 
\label{eq:GenGen}%
\ee
which is a polynomial of degree $|\lambda | + |\mu|$ in $\Q $
and vanishes when $\vQ =1$
$($theorem $\propMainTheorem$~$(\romannumeral1)$ $)$.
}

\proof
\eqref{eq:EigenH} and \eqref{eq:AppCauchy} yield
\ba
Z_{\lambda ,\mu}(\Q ;q,t) 
&=&
\sum_\eta {\Qv^{|\eta|}\over \langle P_\eta , P_\eta \rangle_{q,t} }
P_{\eta}(t^{-\rho};q,t)
\sum_{\nu\leq\lambda \atop \sigma\leq\mu }
\left.
V_{\lambda ,\nu} V_{\mu,\sigma}
H^{\nu^\vee} (x)H^{\sigma^\vee} (x)
P_{\eta}(x;q,t)
\right\vert_{x=t^\rho}
\cr
&=&
\sum_{\nu\leq\lambda \atop \sigma\leq\mu }
V_{\lambda ,\nu} V_{\mu,\sigma}
\left.
H^{\nu^\vee} (x)H^{\sigma^\vee} (x)
\Pi\left(x,\vQ t^{-\rho};q,t\right)
\right\vert_{x=t^\rho}.
\ea
Then, the proposition in subsection \ref{sec:MO}.3
completes the proof. 
\qed

Although
$Z_{\lambda ,\mu}(\Q ;q,t)$
is a power series in $\Q $,
we need its partial sum with degree $|\lambda |+|\mu|$
to calculate ${Z_{\lambda ,\mu}(\Q ;q,t) / Z_{\bullet,\bullet}(\Q ;q,t)}$.

On the assumption that the conjecture 3.4 
is true,
if 
$\vQ =t^N$, $N\in\bZ$ with $0\leq N < \max\{\ell(\lambda ), \ell(\mu)\}$,
then (\ref{eq:GenGen}) 
would vanish 
(theorem $\propMainTheorem$~$(\romannumeral2)$, $|\lambda|+|\mu|\leq 7$ case).

\Section{Towards a resolution of positivity problem }
\label{sec:Positivity}%

We now make a comment on the positivity problem of the GIKV's proposal.
Let $(\bfa; \bq,\bt) := (Q^{-\ha} ;t^{-\ha},-(t/q)^\ha)$.
When $\bfa=\bq^N$ with $N\in\bN$,
the superpolynomial of the homological invariants of the colored Hopf link reduces to 
$
\sum_{i,j \in \bZ} \bq^i \bt^j \dim {\cal H}_{i,j}^{\mathfrak{sl}(N);\lambda,\mu}
$
with certain doubly graded homology ${\cal H}_{i,j}^{\mathfrak{sl}(N);\lambda,\mu}$ \cite{rf:GIKV}.
Therefore, it should be by definition
a polynomial in $\bq $ and $\bt$ with {\it non-negative} integer coefficients.
However, in general, 
$Z^{{\rm inst}}_{\lambda ,\mu}(\vitN ;q,t)\nega{(-1)^{|\lambda |+|\mu|}}$
is not so.
For example, 
\be
-Z^{{\rm inst}}_{3,\bullet}(t^2 ;q,t)
=
 q^3(t^6+t^5)
+q^2(t^5-t^3)
+q  (t^5+t^4).
\ee
A solution to this positivity problem
may be given by replacing the Schur function in \eqref{eq:GIKVZ}
by the Macdonald function $P_\lambda(z;\tildeq,\tildet)$ 
with $\tildet=0$ and appropriately chosen $\tildeq$.
%
%
Let
\be
\tZ_{\lambda ,\mu}(\Q ;q,t) 
:= 
\sum_\eta (-\viQ)^{|\eta|}
P_{\eta^\vee}(q^{\rho};t,q)
P_{\eta}(t^{\rho};q,t)
P_{\lambda}(q^{\eta}t^{\rho};\tildeq,0)
P_{\mu}(q^{\eta}t^{\rho};\tildeq,0),
\ee
and 
$
\tZ^{{\rm inst}}_{\lambda ,\mu}(\Q ;q,t) 
:=
\tZ_{\lambda ,\mu}(\Q ;q,t) 
/
Z_{\bullet,\bullet}(\Q ;q,t). 
$
Note that
$\left.
\tZ^{{\rm inst}}_{\lambda ,\mu}(\Q ;q,t)\right\vert_{\tildeq=0} 
=
Z^{{\rm inst}}_{\lambda ,\mu}(\Q ;q,t) 
$
and
$\tZ^{{\rm inst}}_{1^r ,1^s}(\Q ;q,t) 
=
Z^{{\rm inst}}_{1^r ,1^s}(\Q ;q,t) $,
because of
$P_{\lambda }(x;0,0)=s_{\lambda }(x)$
and 
$P_{1^r}(x;\tildeq,\tildet)=s_{1^r}(x)$,
respectively.
Since
by \eqref{eq:eleUT} and \eqref{eq:MacDef},
\ba
P_{\lambda}(x;\tildeq,\tildet) 
&=& 
\sum_{\mu \leq \lambda }
\tU_{\lambda,\mu} P_\mu (x;q,t),
\qquad
\tU_{\lambda ,\mu}
:=
\sum_{\nu (\lambda \geq \nu \geq\mu )} 
u_{\lambda , \nu }(\tildeq,\tildet) (u^{-1})_{\nu, \mu}(q,t),
\cr
&=&
\sum_{\mu \leq \lambda }
\tV_{\lambda,\mu} e_{\mu^\vee} (x),
\qquad\quad
\tV_{\lambda ,\mu}
:=
\sum_{\nu (\lambda \geq \nu \geq\mu )} 
u_{\lambda , \nu }(\tildeq,\tildet) (a^{-1})_{\nu, \mu},
\ea
all propositions in subsections \ref{sec:GIKV}.3 and \ref{sec:GIKV}.4
hold with these $\tU_{\lambda ,\mu}$ and $\tV_{\lambda ,\mu}$.


If we choose the parameter $\tildeq$ appropriately,
it may overcome the positivity problem.
For example, when $\tildeq=q$,
\be
\left.
-t^{3\over 2}
\tZ^{{\rm inst}}_{3,\bullet}(t^2 ;q,t)
\right\vert_{\tildeq = q}
=
 q^3(t^6+t^5+t^4+t^3)
+q(q+1)(t^5+t^4)
.
\ee
More generally, for $\vQ =t^N$,
\ba
-t^{3\over 2}
\tZ^{{\rm inst}}_{3,\bullet}(\vitN ;q,t)
&=&
\tildeq^3 t^3 {\qbin N3t} 
+\tildeq(\tildeq+1)\left( q\vQ t{\qbin N2t} + t^4{\qbin N3t} \right){\qint 2t}
\cr
&+&
q^3\vQ^2t{\qint Nt} 
+ q^2\vQ t(t^2-1){\qbin N2t}
+q\vQ t^2{\qbin N2t} {\qint 2t} 
+t^6{\qbin N3t}.
\ea
But if we choose $\tildeq=q$, 
then the negative coefficient in the $q^2$-term vanishes as
\ba
\left.
-t^{3\over 2}
\tZ^{{\rm inst}}_{3,\bullet}(\vitN ;q,t)
\right\vert_{\tildeq=q}
&=&
q^3 \left(\vQ^2{\qint Nt} +\vQ {\qbin N2t} {\qint 2t} + t^2{\qbin N3t} \right)t
\cr
&+&
q(q+1)\left( \vQ {\qbin N2t} + t^2{\qbin N3t} \right) t^2{\qint 2t}
+t^6{\qbin N3t}
\cr
&=&
q^3 {\qbin {N+2}3t} +q(q+1)t{\qint 2t}{\qbin {N+1}3t}+ t^3{\qbin N3t}.
\ea
Note that the $q$-integer  
${{\qint Nt} := {1-t^N\over 1-t}}$
and the $q$-binomial coefficient 
${{\qbin Nrt} :=\prod_{i=1}^r{1-t^{N-r+i}\over 1-t^i} }$
are polynomials in $t$ with non-negative integer coefficients for $N,r\in\bN$.
We checked that
$\nega{(-1)^{|\lambda |+s}} \tZ^{{\rm inst}}_{\lambda ,1^s}(\vitN ;q,t)\vert_{\tildeq=q}$
($N\in\bZ_{\geq 0}$) is a polynomial in $q$ and $t$
with non-negative integer coefficients
for $|\lambda |+s\leq 5$ 
(see appendix \appPositivity). 


For non-antisymmetric representations,
the specialization $\tildeq=q$ fails to solve the positivity problem.
For example, 
\be
\left.
\tZ^{{\rm inst}}_{2 ,2}(t^2 ;q,t)
/(qt)^2
\right\vert_{\tildeq=q}
=
q^4t^5 {\qint 2t} + q^3t^3(t^2-1)
+ q^2t^3(t^2+3t+1) + qt^3(t^2+2t+3) + t^2{\qint 3t}.
\ee
However, if we choose
$\tildeq=(1+qt)q/t+p$
with arbitrary $p$, then 
\ba
&&\hskip-48pt
\left.
\tZ^{{\rm inst}}_{2 ,2}(t^2 ;q,t) /(qt)^2
\right\vert_{\tildeq=(1+qt)q/t+p}
=
p^2 t^2 
+ 2p\{q^2t^2 + qt(t^3+t^2+1)+t^3\}
\cr
&+&
q^4t^2(t^4+t^3+1) 
+q^3t(t^4+2t^3+t^2+2)
+ q^2(t^5+t^4+3t^3+t^2+1) 
\cr
&+& qt^2(t^3+2t^2+t+2) + t^2{\qint 3t}
.
\ea
We checked that
$\nega{(-1)^{|\lambda |+|\mu|}}\tZ^{{\rm inst}}_{\lambda ,\mu}(\vitN ;q,t)
\vert_{\tildeq=(1+qt)q/t+p}$
($N\in\bZ_{\geq 0}$) is a polynomial in $q$, $t$ and $p$
with non-negative integer coefficients
for $|\lambda |+|\mu|\leq 5$
(see appendix \appPositivity). 
Although an interpretation from the viewpoint of homological algebra is unclear,
$\tZ^{{\rm inst}}_{\lambda,\mu}(Q ;q,t)\vert_{\tildeq=(1+qt)q/t+p}$
gives a generalization of the colored Hopf link invariants
with a third parameter $p$ in the sense that it reduces to 
$Z^{{\rm inst}}_{\lambda,\mu}(Q ;q,t)$ of GIKV when $\tildeq=0$, i.e. $p=-(1+qt)q/t$.

\section*{Acknowledgments}


We would like to thank J.~Shiraishi, Y.~Yamada and Y.~Yonezawa for discussions. 
This work is partially supported by 
the Grant-in-Aid for Nagoya University Global COE Program, 
`Quest for Fundamental Principles in the Universe: 
from Particles to the Solar System and the Cosmos'. 
This work is also supported in part by Daiko Foundation. 
The work of H.K. is supported in part by 
Grant-in-Aid for Scientific Research [\#19654007] from 
the Japan Ministry of Education, Culture, Sports, Science and Technology.

\Appendix{\appSymm}{Symmetric functions}


Here we recapitulate basic properties of the symmetric functions 
in $x=(x_1,x_2,\cdots)$ \cite{rf:Mac}.
%
%
%
%
%
The monomial symmetric function is defined by
$
m_{\lambda}
:=
\sum_{\sigma}
x_1^{\lambda_{\sigma(1)}}
x_2^{\lambda_{\sigma(2)}} \cdots ,
$
where the summation is over all distinct permutations of 
$(\lambda_1,\lambda_2,\cdots )$.
%
%
The power sum symmetric function is 
$p_n:=\sum_{i=1}^{\infty}x_i^n$.
%
%
The elementary symmetric function defined by 
$$
e_\lambda := e_{\lambda_1}e_{\lambda_2}\cdots
,\qquad
\sum_{r\geq 0} w^r e_r
:=
\prod_{i\geq 1} (1+x_i w)
=
\Exp{-\sum_{n>0}{(-w)^n\over n} p_n}
$$
enjoys
\be
  e_{\lambda^\vee}
  =
  \sum_{\mu\leq\lambda} a_{\lambda\mu} m_{\mu},
  \quad
  a_{\lambda\lambda} = 1,
\quad
  a_{\lambda\mu}\in {\mathbb Z}_{\geq 0}.
\label{eq:eleUT}%
\ee
%
%
For any symmetric functions $f$ and $g$, in power sums $p_n$'s, 
we define a scalar product as
$$
\langle f(p),\, g(p)\rangle_{q,t} 
:= f(p^*)\, g(p)\,
\vert_{{\rm constant\, part}}
p_n^* 
:= 
n {1-q^n \over 1-t^n} 
{\partial \over \partial p_n}.
$$


The Macdonald symmetric function
$P_{\lambda}(x;q,t)$
is uniquely specified by the following orthogonality and normalization:
\ba
  &&
  \langle P_{\lambda}(x;q,t),P_{\mu}(x;q,t)\rangle_{q,t} =0\qquad { \rm if } \;
  \lambda\neq \mu,\\
  &&
  P_{\lambda}(x;q,t)
  =
  \sum_{\mu\leq\lambda} u_{\lambda\mu}(q,t) m_{\mu}(x),
  \quad
  u_{\lambda\lambda}(q,t) = 1,
\quad
  u_{\lambda\mu}(q,t)\in {\mathbb Q}(q,t).
\label{eq:MacDef}%
\ea
Here we used the dominance partial ordering on the Young diagrams defined as
$\lambda\geq\mu \Leftrightarrow |\lambda|=|\mu|$ and
$\lambda_1+\cdots+\lambda_i\geq\mu_1+\cdots+\mu_i$ for all
$i$.
Note that $P(x;q^{-1},t^{-1}) = P(x;q,t)$.
%
%
The scalar product is given by 
$$
\langle P_{\lambda}(x;q,t),P_{\lambda}(x;q,t)\rangle_{q,t} 
=
\prod_{(i,j)\in\lambda}
{
1-q^{\lambda_i-j+1} t^{\lambda_j^\vee-i  }
\over 
1-q^{\lambda_i-j  } t^{\lambda_j^\vee-i+1}
}.
$$
We abbreviate it to $\langle P_{\lambda},P_{\lambda}\rangle_{q,t}$.
%
%
The skew-Macdonald symmetric function 
$P_{\lambda/\mu}(x;q,t)$ is defined by
$$
P_{\lambda/\mu}(x;q,t) 
:=
P_{\mu}^*(x;q,t) \, P_\lambda(x;q,t)
{/\langle P_\mu,P_\mu\rangle_{q,t}}.
$$
Let
$x=(x_1,x_2,\cdots)$ and
$y=(y_1,y_2,\cdots)$
be two sets of variables. Then, we have
\be
\sum_\mu
P_{\lambda/\mu} (x;q,t) 
P_{\mu/\nu} (y;q,t) 
= 
P_{\lambda/\nu} (x,y;q,t),
\label{eq:AppaddSkewMacdonald}
\ee
where
$P_{\lambda/\nu} (x,y;q,t)$
denotes the skew-Macdonald function in the set of variables 
\break
$(x_1,x_2,\cdots,y_1,y_2,\cdots)$.
%
%
%
%
The following Cauchy formula is especially important:
\be
\sum_\lambda
{\langle P_\mu,P_\mu\rangle_{q,t}
\over \langle P_\lambda,P_\lambda\rangle_{q,t}}
P_{\lambda/\mu}(x;q,t) P_\lambda(y;q,t)
=
\Pi(x,y;q,t)P_\mu(y;q,t)
\label{eq:AppCauchy}%
\ee
with
$$
\Pi(x,y;q,t)
:=
\Exp
{
\sum_{n>0}{ 1 \over n}{1-t^n \over 1-q^n} p_n(x) p_n(y)
}.
$$




We denote
\be
p_n(q^\lambda t^\rho) 
:=
\sum_{i=1}^{\ell(\lambda)} (q^{n\lambda_i}-1)t^{n(\ha-i)}
+ {1 \over t^{n\over 2} - t^{-{n\over 2}} }
=
\sum_{i=1}^N q^{n\lambda_i}t^{n(\ha-i)}
+ {t^{-nN} \over t^{n\over 2} - t^{-{n\over 2}} }
\label{eq:powersumI}%
\ee
for any $N\geq \ell(\lambda)$. 
Let
$p_n(x,y) := p_n(x) + p_n(y)$;
then, 
\be
p_n(c q^\lambda t^\rho, L t^{-\rho} ) 
=
c^n\sum_{i=1}^{\ell(\lambda)} (q^{n\lambda_i}-1)t^{n(\ha-i)}
+ {c^n - L^n \over t^{n\over 2} - t^{-{n\over 2}} },
\qquad 
c,L \in\bC. 
\ee
%
%
Note that \cite{rf:AK05}
\be
P_{\mu^\vee}(-t^{\lambda^\vee}q^{\rho};t,q)
=
{v^{|\mu|} \over \langle P_\mu , P_\mu \rangle_{q,t} }
P_{\mu} (q^{-\lambda}t^{-\rho};q,t).
\label{eq:Pdual}%
\ee
The Macdonald function in the power sums
$p_n = (1-L^n)/(t^{{n\over 2}}- t^{-{n\over 2}})$ is
\cite{rf:Mac}(Ch.\ VI.6) 
\be
P_\lambda(t^\rho, L t^{-\rho};q,t) 
=
\prod_{(i,j)\in\lambda}
(-1)t^\ha q^{j-1}
{
1-L q^{1-j} t^{i-1}
\over 
1-q^{\lambda_i-j  } t^{\lambda_j^\vee-i+1}
}
\label{eq:Specialization}%
\ee
for a generic $L\in\bC$.
%
%
Note that
\ba
P_{\lambda }(t^{\rho};q,t)
P_{\lambda ^\vee}(q^{\rho};t,q)
&=&
P_{\lambda }(t^{-\rho};q,t)
P_{\lambda ^\vee}(q^{-\rho};t,q)
\cr
&=&
\prod_{(i,j)\in\lambda}
{
(q/t)^\ha
\over 
(1-q^{-\lambda_i+j  } t^{-\lambda_j^\vee+i-1})
(1-q^{ \lambda_i-j+1} t^{ \lambda_j^\vee-i  })
}.
\label{eq:PPSpecialization}%
\ea


If $L=t^{-N}$ with $N\in\bN$ and $N\geq\ell(\lambda)$, then
$$
p_n(q^\lambda t^\rho, t^{-N-\rho} ) 
= \sum_{i=1}^N q^{n\lambda_i} t^{n(\ha-i)}
$$
is the power sum symmetric polynomial in $N$ variables
$\{q^{\lambda_i}t^{\ha-i}\}_{1\leq i \leq N}$;
hence,
$
P_\lambda(t^\rho, t^{-N-\rho};q,t) 
$ 
reduces to the Macdonald symmetric polynomial in $N$ variables.
Therefore,
$$
P_\lambda(t^\rho, t^{-N-\rho};q,t) =0
\quad{\rm for}\quad
\ell(\lambda) > N\in\bN.
$$
For $N\in\bN$,
there is a symmetry \cite{rf:Mac}(Ch.\ VI.6) 
\be
P_\lambda(      t^\rho,t^{-N-\rho};q,t)
P_\mu (q^\lambda t^\rho,t^{-N-\rho};q,t)
=
P_\mu(      t^\rho,t^{-N-\rho};q,t)
P_\lambda (q^\mu t^\rho,t^{-N-\rho};q,t).
\label{eq:Wsymm}
\ee

The Hall-Littlewood and Schur functions are defined by
$P_{\lambda}(x;t):=P_{\lambda}(x;0,t)$ and
$s_{\lambda}(x):=P_{\lambda}(x;q,q)$, respectively.
Note that $p_n(x)$, $s_\lambda (x)$ and $e_r(x)=P_{1^r}(x;q,t)$
are symmetric functions in $x$ with non-negative integer coefficients.
%
Note that 
$$
\sum_{s=0}^r e_{r-s}(x) e_s(y)=e_r(x,y).
$$
%
For $\lambda = 1^r$, 
\eqref{eq:Specialization} reduces to
\be
e_r(t^\rho,Lt^{-\rho})
=
\prod_{i=1}^r
(-1)t^\ha
{
1-L t^{i-1}
\over 
1- t^i 
}.
\label{eq:SpecializationE}%
\ee 
%
Note also that
$$
P_\lambda (t^\rho;t) 
= \delta_{\lambda ,1^r} 
e_r(t^\rho).
$$
%
%
The $q$-integer  
and the $q$-binomial coefficient 
$$
{\qint Nt} := {1-t^N\over 1-t},
\qquad
{\qbin Nrt}
:=
\prod_{i=1}^r
{
1-t^{N-r+i}
\over 
1-t^i
}
$$
are polynomials in $t$ with non-negative integer coefficients for $N,r\in\bN$.
Note that
$$
t^{-{n\over 2}}p_n(t^{N+\rho},t^{-\rho})
=
{ 1 - t^{nN}\over 1-t^n }
=
{\qint N{t^n}},
\qquad
t^{-{r\over 2}}
e_r (t^{N+\rho},t^{-\rho})
=
t^{{r(r-1)\over 2}}
{\qbin Nrt}.
$$

\Appendix{\appProofHL}{Proof of \eqref{eq:anotherScalarHL} }


Here we prove \eqref{eq:anotherScalarHL}.
%
For $\sigma\in\cS_N$, let
$d(\sigma):=\#\left\{\ (i,j)\ |\ i<j,\ \sigma(i)>\sigma(j)\ \right\}$
be the inversion number and let
\be
\Delta_\sigma 
:=
t^{d(\sigma)}
\Exp
{
\sum_{n>0}{t^n-t^{-n}\over n}
\sum_{i<j \atop \sigma(i) > \sigma(j)}
{x_{\sigma(i)}^n\over x_{\sigma(j)}^n}
}.
\ee
Then,
\be
\Dp(\bar x;t)^{-1} \Delta_\sigma 
= 
t^{d(\sigma)}
\Exp
{
\sum_{n>0}{1-t^{n}\over n}
\left(
\sum_{i<j \atop \sigma(i) < \sigma(j)}
{ x_{\sigma(j)}^n\over x_{\sigma(i)}^n }
-
t^{-n}\sum_{i<j \atop \sigma(i) > \sigma(j)}
{ x_{\sigma(i)}^n\over x_{\sigma(j)}^n }
\right)
},
\ee
which is a formal power series in $\{x_j/x_i\}_{i<j}$
and is equivalent to 
$\sigma \Dp(\bar x;t)^{-1}
=
\prod_{i < j}
{1-tx_{\sigma(j)} /x_{\sigma(i)} \over 1-x_{\sigma(j)} /x_{\sigma(i)} }
$
by the analytic continuation.
Since $P_\lambda (x;t)$ in \eqref{eq:defHL} is a polynomial in $x$,  
it follows that
\ba
P_\lambda (x;t) 
&=&
v_\lambda ^{-1}(t) \Dp(\bar x;t)^{-1}
\sum_{\sigma \in\cS _N}
 \Delta_\sigma  \prod_{i=1}^N x_{\sigma(i)}^{\lambda_i}
\ea
and
\be
v_\lambda (t) P_\lambda (x;t) \Dp(\bar x;t) 
=
\sum_{\nu\geq\lambda }
\sum_{\sigma \in\cS _N}
u_{\lambda ,\nu}^\sigma
\prod_{i=1}^N x_{\sigma(i)}^{\nu_i},
\qquad
u_{\lambda ,\lambda }^\sigma = t^{d(\sigma)}.
\ee
Here 
$\nu$ 
is a sequence of $N$ integers
$\nu = (\nu_1,\nu_2,\cdots,\nu_N)$, $\nu_i\in\bZ$
and $\geq$ is the dominance partial ordering defined as
$\nu\geq\lambda \Leftrightarrow \sum_i\nu_i = \sum_i \lambda_i $ and
$\nu_1 + \cdots + \nu_k \geq \lambda_1+\cdots+\lambda_k$ for all
$k$.
Thus, we have
\ba
v_\lambda (t) 
\langle P_\lambda (x;t) , x^\lambda \rangle_{N;t}''
&=& 1,
\cr
v_\lambda (t) 
\langle P_\lambda (x;t) , \prod_i x_{\sigma(i)}^{\lambda_i} \rangle_{N;t}''
&=& 
t^{d(\sigma )},
\cr
v_\lambda (t) 
\langle P_\lambda (x;t) , \prod_i x_{\sigma(i)}^{\mu_i} \rangle_{N;t}''
&=& 0,
\qquad
\mu<\lambda.
\ea
Therefore, by \eqref{eq:MacDef} we conclude that
\be
v_\lambda (t) 
\langle P_\lambda (x;t) , P_\mu(x;t)\rangle_{N;t}''
= 
\delta_{\lambda, \mu} {\qintf Nt} .
\ee
Here we use the identity
$
\sum_{\sigma\in\cS_N} t^{d(\sigma )} = {\qintf Nt} 
$,
which is proved by the induction in $N$.
This completes the proof of \eqref{eq:anotherScalarHL}.


\Appendix{\appMacOp}{Macdonald operators}


Here we define (higher order) Macdonald operators
which are compatible with tending the number of variables to infinity
\cite{rf:Shi05}\cite{rf:Yam}\cite{rf:AK09}.
For each integer $r$ such that $0\leq r \leq N$,
let $D_N^r$ be the Macdonald operators in $N$ variables 
$x=(x_1,\cdots,x_N)$,
$D_N^0 := 1$ and 
\be
D_N^r := t^{r(r-1)/ 2}
\sum_{I\atop\sharp I = r} \prod_{{i\in I \atop j \notin I}}
{tx_i-x_j\over x_i-x_j}
\prod_{i\in I}{\Diff i},
\label{eq:DNr}%
\ee
summed over all $r$-element subsets $I$ of $\{1,2,\cdots,N\}$.
We set $D_N^r:=0$, $r>N$.
Here
${\Diff{ }}$
is the $q$-shift operator such that
${\Diff{ }} f(x) = f(qx)$. 
Let
$D_N(\tw) := \sum_{r=0}^N D_N^r \tw^r$;
then, the Macdonald polynomial is the eigenfunction for $D_N$:
\ba
&&D_N(\tw) P_\lambda(x;q,t) = P_\lambda(x;q,t) \varepsilon_{N,\lambda},
\cr
&&
\varepsilon_{N,\lambda} := \prod_{i=1}^N (1+\tw q^{\lambda_i}t^{N-i})
= \sum_{r=0}^N \tw^r e_r(q^{\lambda}t^{N-\ha+\rho},t^{\ha-\rho}).
\ea
Therefore, $D_N^r$ are simultaneously diagonalized by the Macdonald 
polynomials
\be
D_N^r P_\lambda(x;q,t) = P_\lambda(x;q,t)  e_r(q^{\lambda}t^{N-\ha+\rho},t^{\ha-\rho});
\ee
thus, $D_N^r$ commute with each other $[D_N^r,D_N^s]=0$ 
on the space of the symmetric function in $N$ variables.
%
%
Note that $D_N^r$ is not compatible with 
the restriction of the variables defined by setting $x_N=0$:
\ba
D_N^r\vert_{x_N=0}
&=&
t^r D_{N-1}^r
+
t^{r-1} D_{N-1}^{r-1},
\cr
D_N(\tw)\vert_{x_N=0}
&=&
(1+\tw) D_{N-1}(t\tw).
\ea
So we need to modify it to take $N\rightarrow\infty$. 
By using 
$
p_n(t^{\rho-\ha})={1/(t^n-1)}
$
and
\break
$
\Exp{ \sum_{n>0}{1\over n}{(-\tw)^n\over 1-t^n} }
=
\sum_{m\geq 0} t^{-{m\over 2}} \tw^m e_m(t^\rho),
$ 
let
\ba
H_N 
&:=&
D_N 
\Exp{ \sum_{n>0}{1\over n}{(-\tw)^n\over 1-t^n} }
=:
\sum_{r\geq 0} \w^r H_N^r,
\qquad
\w:=\tw t^{N-\ha},
\cr
H_N^r &=&
\sum_{s=0}^{\min(r,N)} 
t^{{s\over 2}-rN} e_{r-s}(t^\rho)
D_N^s,
\qquad
r=0,1,2,\cdots.
\label{eq:HNandHNr}%
\ea
Then, by \eqref{eq:powersumI},
we obtain
\ba
E_{N,\lambda} 
:=
\Exp{ \sum_{n>0}{1\over n}{(-\tw)^n\over 1-t^n} }
\varepsilon_{N,\lambda} 
=
\sum_{r\geq 0} \w^r e_r(q^\lambda t^\rho),
\ea
which is independent of $N$ for any $N\geq \ell(\lambda)$. 
Thus,
\ba
H_N P_\lambda(x;q,t) &=& P_\lambda(x;q,t) E_{N,\lambda},
\cr
H_N^r P_\lambda(x;q,t) &=& P_\lambda(x;q,t) e_r(q^\lambda t^\rho).
\ea

\Appendix{\appProofMO}{Proof of \eqref{eq:EigenH} }


Here we prove \eqref{eq:EigenH}
by comparing $H^r$ with 
$H_N^r$ in \eqref{eq:HNandHNr}.
In this subsection, we suppose that the number of variables 
$x=(x_1,\cdots,x_N)$ 
is finite by setting $x_i=0$, $i\geq N+1$ and 
$p_n(x) = \sum_{i=1}^N x_i^n$, $n\in\bN$.
For each integer $r$,
we denote
$\widetilde H_N^r :=t^{rN}H^r/e_r(t^\rho)$,
\ba
\widetilde H_N^r 
&=&
\oint
\prod_{\alpha =1}^r{d\z_\alpha \over 2\pi i\z_\alpha }
\prod_{\alpha =1}^r \prod_{j=1}^N
{t-x_j \z_\alpha \over 1-x_j \z_\alpha}
\prod_{\alpha <\beta }{1-\z_\alpha /\z_\beta \over 1-\z_\alpha /\z_\beta t}
\Exp{
\sum_{n>0}(q^n-1)
\sum_{\alpha =1}^r\z_{\alpha }^{-n}{\partial\over \partial p_n}
}.~~~~~
\ea
Here $\z_\alpha $ and $p_n$ are formal parameters.
For abbreviation,  we write $1/(1-\z)$ instead of $\sum_{n\geq 0}\z^n$
for the formal parameter $\z$.
Note that 
the operators $\widetilde H_N^r$ are compatible with 
the restriction of the variables defined by setting $x_N=0$,
$\widetilde H_{N-1}^r = t^{-r}\widetilde H_N^r\vert_{x_N=0}$.
We also denote 
$\widetilde H_{N-1,(i)}^r := t^{-r} \widetilde H_N^r\vert_{x_i=0}$.

Instead of the rational function
${(tx_i-x_j)/(x_i-x_j)}$,
we use
\ba
\txx{i}{j}
&:=& 
{1-tx_i/x_j\over 1-x_i/x_j}
:=(1-tx_i/x_j)\sum_{n\geq 0}(x_i/x_j)^n,
\cr
\txx{j}{i}
&:=&
{t-x_i/x_j\over 1-x_i/x_j}
:=(t-x_i/x_j)\sum_{n\geq 0}(x_i/x_j)^n,
\qquad  i>j.
\ea
\def\baaiwake{
\be
\txx{i}{j}:= 
\left\{
  \begin{array}{lc}
{\ds{t-x_j/x_i\over 1-x_j/x_i}}
:=(t-x_j/x_i)\sum_{n\geq 0}(x_j/x_i)^n,
\quad & i>j,
\\
{\ds{1-tx_i/x_j\over 1-x_i/x_j}}
:=(1-tx_i/x_j)\sum_{n\geq 0}(x_i/x_j)^n,
\quad & i<j.
    \end{array}
\right.
\ee
}
Then, we have the following recurrence relation.
\\
{\bf Lemma~D.1.}~~{\it
Suppose $x_i\neq x_j$ if $i\neq j$.
For any $N\in\bZ_{\geq 0}$ and $r\in\bN$,
\be
\widetilde H_N^r =
\widetilde H_N^{r-1}
+t^{r-1}(t-1)\sum_{i=1\atop (x_i\neq 0)}^N 
\widetilde H_{N-1,(i)}^{r-1}
\prod_{j(\neq i)} \txx{i}{j}
\Diff{i}.
\label{eq:recurH}%
\ee
}
\proof
For 
$p_n = \sum_i x_i^n$, $n\in\bN$,
since
${\Diff i} p_n = \left((q^n-1) x_i^n + p_n\right)$, 
we have for any function $f$ in $p$
\be
\Diff{i} f(p(x)) = 
\Exp{ \sum_{N>0}(q^n-1)\z^{-n} {\partial\over \partial p_n} }
f(p(x))\vert_{\z=x_i^{-1}}.
\ee
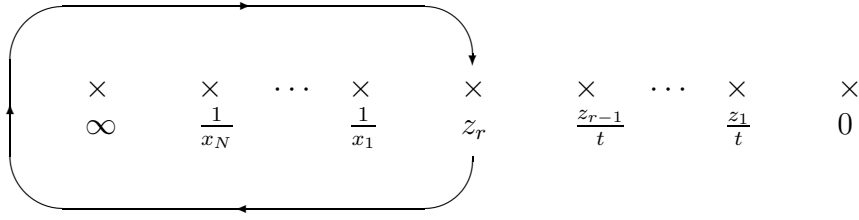
The constant term in $\z_r$ is represented as the contour integral surrounding 
$\infty$
and is written by the summation of the residues at 
$\z_r=\infty$ and $1/x_i$ with $x_i\neq 0$ as follows [Fig.~1]:
\FigContour
\ba
\widetilde H_N^r 
&=&
\oint
\prod_{\alpha =1}^{r-1}{d\z_\alpha \over 2\pi i\z_\alpha }
\prod_{\alpha =1}^{r-1} \prod_{j=1}^N
{t-x_j\z_\alpha \over 1-x_j\z_\alpha }
\prod_{\alpha <\beta }^{r-1}{1-\z_\alpha /\z_\beta \over 1-\z_\alpha /\z_\beta t}
\Exp
{
\sum_{n>0}(q^n-1)
\sum_{\alpha =1}^{r-1}\z_{\alpha }^{-n}
{\partial\over \partial p_n}
}
\cr
&+&
(t-1)\sum_{i=1}^N
\oint
\prod_{\alpha =1}^{r-1}{d\z_\alpha \over 2\pi i\z_\alpha }
\prod_{j(\neq i)}
\txx{i}{j}
\prod_{\alpha =1}^{r-1}{1-\z_\alpha x_i \over 1-\z_\alpha x_i/t}
\Diff{i}
\cr
&\times&
\prod_{\alpha =1}^{r-1} \prod_{j=1}^N
{t-x_j \z_\alpha \over 1-x_j \z_\alpha}
\prod_{\alpha <\beta }^{r-1}{1-\z_\alpha /\z_\beta \over 1-\z_\alpha /\z_\beta t}
\Exp
{
\sum_{n>0}(q^n-1)
\sum_{\alpha =1}^{r-1}\z_{\alpha }^{-n}
{\partial\over \partial p_n}
}
\cr
&=&
\widetilde H^{r-1}
+t^{r-1}(t-1) 
\sum_{i=1}^N
\oint
\prod_{\alpha =1}^{r-1}{d\z_\alpha \over 2\pi i\z_\alpha }
\prod_{j(\neq i)}
\txx{i}{j}\cdot
\Diff{i}
\cr
&&\hskip30pt
\times
\prod_{\alpha =1}^{r-1} \prod_{j(\neq i)}
{t-x_j \z_\alpha \over 1-x_j \z_\alpha}
\prod_{\alpha <\beta }^{r-1}{1-\z_\alpha /\z_\beta \over 1-\z_\alpha /\z_\beta t}
\Exp{
\sum_{n>0}(q^n-1)
\sum_{\alpha =1}^{r-1}\z_{\alpha }^{-n}
{\partial\over \partial p_n}
}.
\ea
\qed

Let us denote
$\Ii{i}\oplus \{i\} \oplus \Ji{i}:= \{1,2,\cdots,N\}$. 
Let
\be
D_{N-1,(k)}^r :=
t^{r(r-1)/ 2}
\sum_{\Ii{k} \atop \sharp \Ii{k} = r}
\prod_{{i\in \Ii{k} \atop j \in \Ji{k}}} 
\txx{i}{j}
\prod_{i\in \Ii{k}}{\Diff i};
\ee
then, we have the following recurrence relation.
\\
{\bf Lemma~D.2.}~~{\it
For any $r=0,1,\cdots,N-1$,
\be
\sum_{i=1}^N 
\prod_{j(\neq i)} \txx{i}{j}
\cdot
D_{N-1,(i)}^r
\Diff{i}
=
t^{-r} {t^{r+1} -1\over t-1}
D_N^{r+1}.
\ee
}
\proof
\ba
t^{-{r(r-1)/2}}\times
\lhs
&=&
\sum_{i=1}^N 
\prod_{j(\neq i)} \txx{i}{j}
\cdot
\Diff{i}
\sum_{\Ii{i}\atop \sharp \Ii{i}=r} 
\prod_{k\in \Ii{i} \atop \ell\in \Ji{i}} 
\txx{k}{\ell}
\prod_{k\in \Ii{i}}
\Diff{k}
\cr
&=&
\sum_{i=1}^N 
\sum_{\Ii{i}}
\prod_{j\in \Ii{i}\cup\Ji{i}} 
\txx{i}{j}
\cdot
\Diff{i}
\prod_{k\in \Ii{i} \atop \ell\in \Ji{i}} 
\txx{k}{\ell}
\prod_{k\in \Ii{i}}
\Diff{k}
\cr
&=&
\sum_{i=1}^N 
\sum_{\Ii{i}}
\prod_{j\in \Ii{i}} 
\txx{i}{j}
\prod_{k\in I \atop \ell\in \Ji{i}} 
\txx{k}{\ell}
\prod_{k\in I}
\Diff{k},
\qquad I:=\Ii{i}\oplus \{i\},
\cr
&=& 
\sum_{I\atop\sharp I = r+1} 
\sum_{i\in I}
\prod_{j\in I\atop j\neq i} \txx{i}{j}
\prod_{k\in I \atop \ell\notin I} 
\txx{k}{\ell}
\prod_{k\in I}
\Diff{k}.
\ea
Thus, it is sufficient to show that
$
\sum_{i=1}^{r+1}
\prod_{j(\neq i)} \txx{i}{j}
=
\sum_{i=0}^{r} t^i,
$ 
which is proved as follows.

\noindent (i)
Since the residues at $x_i=x_j$ vanish, so the \lhs\ 
is a constant.

\noindent (ii)
By putting $x_i = \epsilon^i$ and taking $\epsilon\rightarrow \infty$, 
we obtain the \rhs.
\qed

Hence, we have
\\
{\bf Lemma~D.3.}~~{\it
\be
\widetilde H_N^r =
\sum_{k=0}^r D_N^k
\prod_{i=0}^{k-1} (t^{r-i}-1).
\label{eq:HbyD}%
\ee
}
\proof
We proceed by induction on $r$.
When $r=0$, since $\widetilde H_N^0=1$, 
(\ref{eq:HbyD}) holds.
So assume that the result is true for $r-1\geq 0$.
{}From (\ref{eq:recurH}) we obtain
\ba
\widetilde H_N^r 
&=&
\sum_{k=0}^{r-1} 
\prod_{s=0}^{k-1} (t^{r-s-1}-1)
\left(
D_N^k
+
t^{r-1} (t-1)
\sum_{i=1}^N D_{N-1,(i)}^k
\prod_{j(\neq i)} \txx{i}{j}
\Diff{i}
\right)
\cr
&=&
\sum_{k=0}^{r-1} 
\prod_{s=0}^{k-1} (t^{r-s-1}-1)
\left(
D_N^k
+
(t^r-t^{r-k-1}) 
D_N^{k+1}
\right)
\cr
&=&
\left(
\sum_{k=0}^{r-1} 
(t^{r-k}-1) 
+ 
\sum_{k=1}^{r} 
(t^r-t^{r-k})
\right)
D_N^k
\prod_{s=1}^{k-1} (t^{r-s}-1).
\ea
\qed

Therefore,
\\
{\bf Proposition~D.4.}~~{\it
\ba
H^r &=&
\sum_{s=0}^{\min(r,N)} 
t^{{s\over 2}-rN} 
e_{r-s}(t^\rho)
D_N^s,
\qquad
r=0,1,2,\cdots.
\ea
}
This completes the proof of (\ref{eq:EigenH})
by taking the limit $N\rightarrow\infty$.


Since the Macdonald functions for all partitions 
form a basis of the ring of symmetric functions,
$H^r$ commute with each other on the space of symmetric functions.


\Appendix{\appTorus}{Torus knot}


Not only for the Hopf link but the homological invariants for other link 
may be related with the refined topological vertex or the Macdonald polynomial. 
For example, in \cite{rf:DGR},
a reduced polynomial for the torus knot $T_{m,n}$ is conjectured as 
\be
{\cal P}(T_n)
:=
\lim_{m\rightarrow \infty} {\cal P}(T_{m,n})
=
{ 1-a^2 \bt u \over 1- \bt^{-2} u^2}
{ 1-a^2 \bt u^2 \over 1- \bt^{-2} u^3}
\cdots
{ 1-a^2 \bt u^{n-1} \over 1- \bt^{-2} u^n}
\ee
with 
$u:=(\bq \bt)^2$.
But this is nothing but the following specialization of the Macdonald polynomial: 
\be
P_{(n-1)}(Qt^{\ha+\rho},t^{-\ha-\rho};q,t)
=
\prod_{i=0}^{n-2}
{ 1-q^i t Q \over 1-q^i t}
\ee
with
$(Q;q,t)=(a^2\bt/\bq^2;\bq^2 \bt^2,\bq^4 \bt^2)$.
Note that
$p_n (Qt^{\ha+\rho},t^{-\ha-\rho} ) = {(1-(Qt)^n)/ (1-t^n) }$.


\Appendix{\appPositivity}{Example for $\tZ^{{\rm inst}}_{\lambda ,\mu}(t^N ;q,t)$ }


Here we list the example for 
$\tZ^{{\rm inst}}_{\lambda, \mu}(t^N ;q,t)$
in section \ref{sec:Positivity}
with $Q=t^N$,
which is written as a linear combination of 
the $q$-binomial coefficients ${\qbin Nrt}$
with $r\leq |\lambda|+|\mu|$.
When $\tildeq=(1+qt)q/t+p$ for $|\lambda|+|\mu|\leq 5$,
we checked that any coefficient of them is
a polynomial in $q$, $t$ and $p$
with non-negative integer coefficients, 
so is $\tZ^{{\rm inst}}_{\lambda, \mu}(t^N ;q,t)$.
The same is true 
when $\tildeq=q$,
i.e. $p=q-(1+qt)q/t$,
for $(\lambda,\mu)=(\lambda,1^s)$ with $|\lambda|+s\leq 5$.

\be 
-t^\ha \tZ^{{\rm inst}}_{1 ,\bullet}\tNqt
=
t{\qint{N}t}
\ee
\be 
\tZ^{{\rm inst}}_{1 , 1}\tNqt
=
(qQ+{\qint{N}t}){\qint{N}t}
\ee
\be 
t
\begin{pmatrix}
\tZ^{{\rm inst}}_{2 ,\bullet}\tNqt \cr
\tZ^{{\rm inst}}_{1^2 ,\bullet}\tNqt
\end{pmatrix}
=
\begin{pmatrix}
~~1~~ & ~~\tildeq +t~~ \cr
0 & 1
\end{pmatrix}
\begin{pmatrix}
qQt{\qint Nt}
\cr
t^2{\qbin N2t}
\end{pmatrix},
\ee
\be 
-t^\ha
\begin{pmatrix}
\tZ^{{\rm inst}}_{2 ,1}\tNqt \cr
\tZ^{{\rm inst}}_{1^2 ,1}\tNqt 
\end{pmatrix}
=
\begin{pmatrix}
~~q~~ & \tildeq+q(t-1)+{\qint 2t} & ~~\tildeq+t~~ \cr
0 & 1 & 1
\end{pmatrix}
\begin{pmatrix}
q^2Q^2
{\qint Nt}
\cr
qQ
{\qint 2t}
{\qbin N2t}
\cr
t^2
{\qint 3t}
{\qbin N3t}
\end{pmatrix},
\ee
\be 
\begin{pmatrix}
\tZ^{{\rm inst}}_{2 ,1^2}\tNqt \cr
\tZ^{{\rm inst}}_{1^2 ,1^2}\tNqt 
\end{pmatrix}
=
\begin{pmatrix}
~~\tildeq+qt{\qint 2t}+t~~ & (\tildeq+q(t-1)){\qint 2t}+{\qint 3t} & ~~\tildeq+t~~ \cr
1 & 1 & 1 \cr
\end{pmatrix}
\begin{pmatrix}
{q^2Q^2\over t^2} 
{\qbin N2t}
\cr
{qQ\over t} 
{\qint 3t}
{\qbin N3t}
\cr
t^2
{\qbin 42t}
{\qbin N4t}
\end{pmatrix},
\ee
\be 
-t^{-\ha}
\begin{pmatrix}
\tZ^{{\rm inst}}_{2 ,1^3}\tNqt \cr
\tZ^{{\rm inst}}_{1^2 ,1^3}\tNqt 
\end{pmatrix}
=
\begin{pmatrix}
~~\tildeq+(qt+1)t~~ & (\tildeq+q(t-1)){\qint 3t}+{\qint 4t} & ~~\tildeq+t~~ \cr
1 & {\qint 3t} & 1 \cr
\end{pmatrix}
\begin{pmatrix}
{q^2Q^2\over t^4} 
{\qint 3t}
{\qbin N3t}
\cr
{qQ\over t^2} 
{\qint 4t}
{\qbin N4t}
\cr
t^2
{\qbin 53t}
{\qbin N5t}
\end{pmatrix}.
\ee
Note that ${\qint 4t}=(t^2+1){\qint 2t}$.
Since
$\{\tildeq+q(t-1)\}{\qint 2t}$ 
with
$\tildeq=q+p$ or $(1+qt)q/t+p$
is a polynomial in $q$, $t$ and $p$
with positive integer coefficients,
so is 
$\tZ^{{\rm inst}}_{\lambda ,1^s}\tNqt $
for $|\lambda|=2$ and $s=0,1,2,3$.

\be 
-
t^{{3\over 2}}
\begin{pmatrix}
\tZ^{{\rm inst}}_{3 ,\bullet}\tNqt \cr
\tZ^{{\rm inst}}_{(2,1) ,\bullet}\tNqt \cr
\tZ^{{\rm inst}}_{1^3 ,\bullet}\tNqt 
\end{pmatrix}
=
M_{3,0}
\begin{pmatrix}
q^2Q^2
{\qint Nt} \cr
qQt
{\qint 2t} 
{\qbin N2t} \cr
t^3
{\qbin N3t} 
\end{pmatrix},
\qquad
t
\begin{pmatrix}
\tZ^{{\rm inst}}_{3 ,1}\tNqt \cr
\tZ^{{\rm inst}}_{(2,1) ,1}\tNqt \cr
\tZ^{{\rm inst}}_{1^3 ,1}\tNqt 
\end{pmatrix}
=
M_{3,1}
\begin{pmatrix}
q^4Q^3
{\qint Nt} \cr
{q^2Q^2 \over t} 
{\qint 2t} 
{\qbin N2t} \cr
qQ
{\qint 3t} 
{\qbin N3t} \cr
t^{3}
{\qint 4t} 
{\qbin N4t} 
\end{pmatrix},
\ee
\be
-t^\ha
\begin{pmatrix}
\tZ^{{\rm inst}}_{3 ,1^2}\tNqt \cr
\tZ^{{\rm inst}}_{(2,1) ,1^2}\tNqt \cr
\tZ^{{\rm inst}}_{1^3 ,1^2}\tNqt 
\end{pmatrix}
=
M_{3,2}
\begin{pmatrix}
q^4Q^3
{\qint 2t} 
{\qbin N2t} \cr
{q^2Q^2 \over t^3}
{\qint 3t} 
{\qbin N3t} \cr
{qQ\over t}
{\qint 3t}{\qint 4t} 
{\qbin N4t} \cr
t^3
{\qbin 52t} 
{\qbin N5t} 
\end{pmatrix},
\qquad
t^2
\begin{pmatrix}
\tZ^{{\rm inst}}_{4 ,\bullet}\tNqt \cr
\tZ^{{\rm inst}}_{(3,1) ,\bullet}\tNqt \cr
\tZ^{{\rm inst}}_{2^2 ,\bullet}\tNqt \cr
\tZ^{{\rm inst}}_{(2,1^2) ,\bullet}\tNqt \cr
\tZ^{{\rm inst}}_{1^4 ,\bullet}\tNqt 
\end{pmatrix}
=
M_{4,0}
\begin{pmatrix}
q^3Q^3 
{\qint Nt} \cr
q^2Q^2
{\qbin N2t} \cr
 q Q t {\qint 3t}
{\qbin N3t} \cr
t^4
{\qbin N4t} 
\end{pmatrix},
\ee
with $M_{3,s}$ and $M_{4,0}$ in the next page.
%
Since
\ba
&&
\{\tildeq+q(t^2-1)\}{\qint 3t},
\qquad
\tildeq(\tildeq+1)(qt+1)+q^2(t-1)(qt+{\qint 2t}), 
\cr
&&
\{\tildeq^3 +(\tildeq(\tildeq+1)+q(t-1))(q(t-1)+{\qint 2t}){\qint 2t} \}{\qint 3t},
\ea
with
$\tildeq=q+p$ or $(1+qt)q/t+p$
are polynomials in $q$, $t$ and $p$
with positive integer coefficients,
so is 
$\tZ^{{\rm inst}}_{\lambda, 1^s}\tNqt $
for $|\lambda|=3$ and $s=0,1,2$.
The same is true for $|\lambda|=4,5$ and $|\lambda|+s\leq 5$
with $\tildeq=q$ or $(1+qt)q/t+p$.

\ba
&&\hskip-12pt
t^2\tZ^{{\rm inst}}_{2,2}(\vitN ;q,t) 
=
q^6 Q^3 t{\qint Nt}
+(\tildeq+t)^2 t^4(t^2\!+\!1){\qint 3t}{\qbin N4t} 
\cr
&+&
q^2 Q^2
\left\{
\left\{
2\tildeq t + q(qt+{\qint 2t})(t\!-\!1)
\right\}
q{\qint 2t}
+(\tildeq^2+t)^2+(qt{\qint 2t}+1){\qint 2t}
\right\}
{\qbin N2t} 
\cr
&+&
qQt
 \left\{
( \tildeq+q(t\!-\!1) )( \tildeq+q(t\!-\!1)+{\qint 2t} ){\qint 2t}
+\tildeq^2(t^2\!+\!1) + t({\qint 3t}+1)
\right\}
{\qint 3t}{\qbin N3t} 
.
~~~~~~~~
\ea
Since
$\{\tildeq+q(t-1)\}{\qint 2t}$,
$\{\tildeq+q(t-1)\}{\qint 3t}$,
$\tildeq t +(t-1)q(qt+{\qint 2t})$,
with
$\tildeq=(1+qt)q/t+p$
are polynomials in $q$, $t$ and $p$
with positive integer coefficients,
so is 
$\tZ^{{\rm inst}}_{2,2}\tNqt $.
The same is true for $(\lambda,\mu)=(3,2)$ and $((2,1),2)$
with $(1+qt)q/t+p$.


\begin{landscape}


\be
{}^t M_{3,0}=
\begin{pmatrix}
qt & 0 &~ 0 ~~\cr
\tildeq(\tildeq+1)+q(t-1)+t & 1 &~ 0 ~~\cr
~~ \tildeq^3 +\tildeq(\tildeq+1)t{\qint 2t}+t^3 
~&~ \tildeq+t{\qint 2t} 
~&~ 1 
~~\cr
\end{pmatrix}.
\nonumber
\ee
\be
{}^t M_{3,1}=
\begin{pmatrix}
q^2 
& 0 
&~ 0 
~~\cr
\tildeq(\tildeq+1)(qt+1)+q^2(t-1)(qt+{\qint 2t})+(q{\qint 2t}+1)t 
& qt+1 
&~ 0  
~~\cr
~~ \tildeq^3 +(\tildeq(\tildeq+1)+q(t-1))(q(t-1)+{\qint 2t}){\qint 2t}+t{\qint 3t} 
~&~ \tildeq+(q(t-1)+{\qint 2t}){\qint 2t} 
~&~ 1 ~~\cr
~~ \tildeq^3 +\tildeq(\tildeq+1)t{\qint 2t}+t^3 
~&~ \tildeq+t{\qint 2t} 
~&~ 1 
~~\cr
\end{pmatrix}.
\nonumber
\ee
\ba
&&
{}^t M_{3,2}=
\cr
&&
\begin{pmatrix}
\tildeq(\tildeq+1)+q^2t^2+q(t-1)+t
& t^{-3} 
&~ 0 
~~\cr
~~
\tildeq^3\!+\!\tildeq(\tildeq\!+\!1)(q(t^2\!+\!t\!-\!1)\!+\!{\qint 2t}){\qint 2t}
\!+\!qt(t\!-\!1)(q{\qint 2t}(qt\!+\!{\qint 2t})\!+\!1)
\!+\!qt^3(t\!+\!2)\!+\!t{\qint 3t} 
~&~ \tildeq\!+\!q(t^2\!+\!t\!-\!1){\qint 2t}\!+\!{\qint 2t}^2 
~&~ 1  
~~\cr
\tildeq^3 +(\tildeq(\tildeq+1)+q(t-1))(q(t^2-1)+{\qint 3t})+t(t^2+1)
& \tildeq+q(t^2-1)+{\qint 3t}
&~ 1 
~~\cr
\tildeq^3 +\tildeq(\tildeq+1)t{\qint 2t}+t^3 
& \tildeq+t{\qint 2t} 
&~ 1 
~~\cr
\end{pmatrix}.
\nonumber
\ea

{\footnotesize
\vbox{%
\ba
&&
{}^t M_{4,0}=
\cr
&&
\begin{pmatrix}
q^3t & 0 & 0 &0 &0 \cr
\tildeq^4+\tildeq^2 +\tildeq {\qint 3\tildeq} (q{\qint 2t}+1)t
+q^2 (t-1) {\qint 2t} (qt+{\qint 2t})+(q{\qint 2t}+1)t^2
& \tildeq+{\qint 2t} t q+t
& 1  
& 0 &0\cr
q^2t( \tildeq {\qint 3\tildeq}+(\tildeq^2+1)t )
+( (t-1) q+t ) ( \tildeq {\qint 3\tildeq}+(t-1) q ) {\qint 2t}+t^3
& \tildeq^2+{\qint 2t} (\tildeq+(t-1) q +t)
& \tildeq+t
& 1 &0
\cr
\tildeq^6+\tildeq {\qint 3\tildeq} (\tildeq^2+t^2) t{\qint 3t} 
+(\tildeq^4+\tildeq^2)(t^4+t^2) 
+t^6
&\hskip-12pt
\tildeq^3+ \tildeq(\tildeq+1) t {\qint 3t} +\tildeq t^2 (t^2+1)+ t^3{\qint 3t}
~&\tildeq^2+ \tildeq t {\qint 3t}+ t^4+t^2 
~&\tildeq+t {\qint 3t} 
~&1
\cr
\end{pmatrix}.
\nonumber
\ea
}}

\newpage
\ba
&&\hskip-12pt
t^{5\over 2}\tZ^{{\rm inst}}_{(2,1),2}(\vitN ;q,t) 
=
(\tildeq+t{\qint 2t})(\tildeq+t)(t^2\!+\!1)t^5{\qint 5t}{\qbin N5t} 
+
{q^4Q^3\over t}
 \left\{
\tildeq+q(t\!-\!1)+q^2t^2+ {\qint 2t}
\right\}
{\qint 2t}{\qbin N2t} 
\cr
&+&
{q^2 Q^2\over t}
\left\{
\tildeq^2+
q\left\{
\tildeq + q^2t^2+(qt\!+\!1){\qint 2t}
\right\}
(t\!-\!1)
+\tildeq
\left\{
qt^2(t\!+\!2)+{\qint 3t}+2t
\right\}
+
({\qint 3t}+qt^2{\qint 2t}){\qint 2t}
\right\}
{\qint 3t}{\qbin N3t}
\cr
&+&
q Q
\left\{
\left\{
\tildeq^2 + q(q{\qint 2t}(t\!-\!1)+\tildeq(t\!+\!2)+{\qint 2t}^2)(t\!-\!1)
\right\}
{\qint 3t}
+\tildeq({\qint 3t}^2+{\qint 4t})+t({\qint 3t}^2+{\qint 2t}) 
\right\}
t(t^2\!+\!1){\qint 2t}{\qbin N4t} 
.
\nonumber
\ea

\ba
&&\hskip-12pt
t^{5\over 2}\tZ^{{\rm inst}}_{3,2}(\vitN ;q,t) 
=
 q^{10} Q^4 {\qbin N1t}  t
+
\left\{
 \tildeq^4
 +\tildeq^3 t (t+2)
 +\tildeq^2 t {\qint 2t}^2
 +\tildeq t^2 (2 t+1)
 +t^4
\right\}
  (t^2+1)t^5 {\qint 5t}  {\qbin N5t} 
\cr
&+& 
 {q^4 Q^3\over t} 
 \left\{
 \tildeq^3
+q
\left\{
 \tildeq (\tildeq+2)
 +((q^3 t+1) (q t+{\qint 2t})+q^2 t{\qint 2t} ) 
 \right\}
  (t\!-\!1)
\right.\cr&&~~~~~~~~\left.
+\tildeq^2(q^2 t^2 +t+2) 
 +\tildeq(2 q^2 t^2 +2 t+1) 
 +(q^2 t^2{\qint 2t} +{\qint 3t})
\right\}
{\qint 2t} {\qbin N2t}
\cr
&+&
 {q^2 Q^2\over t} 
 \left\{
\tildeq^4
+ q 
 \left\{
q t (q t+{\qint 2t})+{\qint 2t}
\right\}
(\tildeq^2+q^2(t\!-\!1))  (t\!-\!1) 
{\qint 2t} 
\right.\cr&+&\left.
 q 
 \left\{
 \tildeq^2(\tildeq+1) {\qint 2t}
 +\tildeq((2 q t (q t+ {\qint 2t})+1)
  {\qint 2t}+2 t) 
 +((2 q^2 t^2+q{\qint 2t}^2 +2) t {\qint 2t}+1)
\right\}
 (t\!-\!1)
\right.\cr&+&\left.
  \left\{
   \tildeq^3(q t^2 (t\!+\!2)\!+\!{\qint 3t}\!+\!2 t) 
 \!+\!\tildeq^2 ({\qint 3t}\!+\!{\qint 2t}\!+\!q t^2 (t\!+\!2)) {\qint 2t}
 \!+\!\tildeq((2 t\!+\!1) {\qint 3t}\!+\!qt^3(2 t\!+\!5) ) 
 \!+\!(q t^2{\qint 2t}^2 \!+\!{\qint 4t}\!+\!t) t
\right\}
\right\}
{\qint 3t} {\qbin N3t}
 \cr
&+&
qQ
 \left\{
q^3 {\qint 2t} {\qint 3t}  (t\!-\!1)^3
 \!+\! q^2(\tildeq (\tildeq\!+\!2)\!+\!{\qint 2t}) {\qint 2t} {\qint 3t} (t\!-\!1)^2
 \!+\!q\left\{
  \tildeq^3(t\!+\!2)
 \!+\!\tildeq (2 {\qint 3t}\!+\!t) 
 \!+\!\tildeq^2{\qint 2t} (t\!+\!2) 
 \!+\!t ({\qint 3t}\!+\!1)
\right\}
 {\qint 3t} (t\!-\!1)
 \right.\cr&&~~~~\left.
+\tildeq^4 {\qint 3t}
 +\tildeq^3({\qint 3t}^2+{\qint 4t}) 
 + \tildeq ({\qint 4t} {\qint 2t}+{\qint 3t} (t^2+2)) t
 +({\qint 4t}+t) {\qint 2t}^2 \tildeq^2
 + ({\qint 5t}+t{\qint 2t} )t^2
\right\}
 t {\qint 4t} {\qbin N4t}
.
\nonumber
\ea


\end{landscape}





\end{document}